\documentclass[11pt]{amsart}
\usepackage{amsmath,amssymb}
\newtheorem{theorem}{Theorem}[section]
\newtheorem{proposition}[theorem]{Proposition}
\newtheorem{corollary}[theorem]{Corollary}
\newtheorem{lemma}[theorem]{Lemma}
\theoremstyle{definition}
\newtheorem{definition}[theorem]{Definition}
\newtheorem{assumption}[theorem]{Assumption}

\begin{document}

\title[Mean curvature flow with surgery]{Mean curvature flow with surgery of mean convex surfaces in $\mathbb{R}^3$}
\author{Simon Brendle and Gerhard Huisken}
\address{Department of Mathematics \\ Stanford University \\ Stanford, CA 94305}
\address{Mathematisches Institut \\ Universit\"at T\"ubingen \\ 72076 T\"ubingen \\ Germany}
\begin{abstract}
We define a notion of mean curvature flow with surgery for two-dimensional surfaces in $\mathbb{R}^3$ with positive mean curvature. Our construction relies on the earlier work of Huisken and Sinestrari in the higher dimensional case. One of the main ingredients in the proof is a new estimate for the inscribed radius established by the first author \cite{Brendle2}.
\end{abstract}
\maketitle

\section{Introduction}

The formation of singularities in geometric flows is a central problem in geometric analysis. In the 1990s, Hamilton \cite{Hamilton1} started a program aimed at understanding the singularities of the Ricci flow in dimensions $3$ and $4$. In particular, for Ricci flow on four-manifolds with positive isotropic curvature, Hamilton \cite{Hamilton2} showed that the flow can be extended beyond singularities by means of surgery procedure. In 2002, Perelman \cite{Perelman1}, \cite{Perelman2}, \cite{Perelman3} successfully carried out a surgery construction for the Ricci flow in dimension $3$, and used it to prove the Poincar\'e and Geometrization Conjectures.

It this paper, we focus on the mean curvature flow. In \cite{Huisken-Sinestrari3}, Huisken and Sinestrari defined a notion of mean curvature flow with surgery for two-convex hypersurfaces in $\mathbb{R}^{n+1}$, where $n \geq 3$. We assume that the reader is familiar with that paper. Our goal in this paper is to extend the construction in \cite{Huisken-Sinestrari3} to the case $n=2$:

\begin{theorem}
\label{main.theorem}
Let $M_0$ be a closed, embedded surface in $\mathbb{R}^3$ with positive mean curvature. Then there exists a mean curvature flow with surgeries starting from $M_0$ which terminates after finitely many steps.
\end{theorem}

As in \cite{Huisken-Sinestrari3}, the surgery construction involves three curvature thresholds $H_1,H_2,H_3$, where $H_3 = 10H_2 \gg H_1$. The basic idea is that we let the flow evolve smoothly until the maximum curvature reaches the threshold $H_3$. When that happens, we perform surgeries on necks with a curvature scale comparable to $H_1$. As a result of that, the maximum curvature drops below $H_2$ right after surgery. We then let the flow evolve smoothly until the maximum curvature reaches the threshold $H_3$ again, and repeat the process until the flow becomes extinct. 

If we send the curvature thresholds $H_1,H_2,H_3$ to infinity, the flow with surgery will converge to the level set solution; this follows from results of Head \cite{Head} and Lauer \cite{Lauer}.

Our argument broadly follows the one in \cite{Huisken-Sinestrari3}. However, there are several major differences. One important difference is that the cylindrical estimate in Section 5 of \cite{Huisken-Sinestrari3} fails for $n=2$. To replace the cylindrical estimate, we use an estimate for the inscribed radius established in \cite{Brendle2} (see also \cite{Brendle3} for a recent survey). Given an embedded oriented surface $M$ in $\mathbb{R}^3$ and a point $p \in M$, the inscribed radius at $p$ is defined as the radius of the largest open ball in $\mathbb{R}^3$ which is disjoint from $M$ and touches $M$ at $p$ from the inside. Similarly, the outer radius at $p$ is defined as the radius of the largest open ball in $\mathbb{R}^3$ which is disjoint from $M$ and touches $M$ at $p$ from the outside. Following Sheng and Wang \cite{Sheng-Wang}, an embedded mean convex surface $M$ will be called $\alpha$-noncollapsed if the inscribed radius at each point $p \in M$ is bounded from below by $\frac{\alpha}{H}$, where $H$ denotes the mean curvature at the point $p$. It follows from general results of Brian White \cite{White2},\cite{White3} that every embedded solution of the mean curvature flow with positive mean curvature is $\alpha$-noncollapsed for some uniform constant $\alpha>0$ which is independent of $t$. An alternative proof of that fact was given by Sheng and Wang \cite{Sheng-Wang}. Andrews recently showed that the $\alpha$-noncollapsing condition is preserved by the flow; this uses a maximum principle argument similar in spirit to \cite{Huisken3}. 

The main theorem in \cite{Brendle2} asserts that, for any smooth solution of the mean curvature flow with positive mean curvature, we have a pointwise estimate of the form $\mu \leq (1+\delta) \, H + C(\delta)$. Here, $\mu$ denotes the reciprocal of the inscribed radius, $\delta$ is a given positive number, and $C(\delta)$ is a constant that depends on $\delta$ and the initial data. We show that this estimate still holds in the presence of surgeries, at least for a suitable choice of surgery parameters. This is a subtle issue, as the ratio $\frac{\mu}{H}$ might deteriorate slightly under surgery. To overcome this obstacle, we show that the ratio $\frac{\mu}{H}$ improves immediately prior to surgery. By a suitable choice of the surgery parameters, we can ensure that this improvement in the noncollapsing constant prior to surgery is strong enough to absorb the error terms that arise during each surgery procedure.

Another problem is that the proof of the gradient estimate in Section 6 of \cite{Huisken-Sinestrari3} does not directly carry over to the case $n=2$. To get around this issue, we use a new interior gradient estimate due to Haslhofer and Kleiner \cite{Haslhofer-Kleiner}. The estimate of Haslhofer and Kleiner implies that $|\nabla A| \leq C \, H^2$, provided that the flow is $\alpha$-noncollapsed and has evolved smoothly for a long enough time (cf. Theorem \ref{interior.derivative.estimate} below). On the other hand, we can use the pseudolocality principle to control $|\nabla A|$ shortly after surgery (cf. Proposition \ref{consequence.of.pseudolocality}). By combining these two results, we obtain an estimate for $|\nabla A|$ which is valid at all points in space-time, even in the presence of surgeries (see Proposition \ref{gradient.estimate}).

In Section \ref{overview}, we state a number of auxiliary results. In Section \ref{neck.contin}, we use these auxiliary results to establish an analogue of the crucial Neck Continuation Theorem in \cite{Huisken-Sinestrari3}. We then implement the surgery construction from \cite{Huisken-Sinestrari3}, and complete the proof of Theorem \ref{main.theorem}. Finally, in Sections \ref{proof.of.pseudoloc} -- \ref{proof.of.7.19} we give the proofs of the auxiliary results stated in Section \ref{overview}.

Finally, let us mention some related results. Brian White has obtained several breakthroughs in the analysis of the singularities of mean convex mean curvature flow; see \cite{White1}, \cite{White2}, \cite{White3}, \cite{White4}, \cite{White5}, \cite{White6}. A different approach to the singularity analysis for mean curvature flow in the two-dimensional case was suggested by Colding and Kleiner in \cite{Colding-Kleiner}. Moreover, Wang \cite{Wang} has obtained a classification of translating solutions to the mean curvature flow in dimension $2$. These solutions arise as models for Type II singularities. Finally, the first author has recently obtained a classification of self-similar solutions to the Ricci flow in dimension $3$ under a noncollapsing assumption (see \cite{Brendle1}).

We are grateful Brian White for discussions concerning the pseudolocality property for the mean curvature flow. We thank the referees for their careful reading of the original manuscript, and for valuable comments.

\section{Overview of some auxiliary results}

\label{overview}

In this section, we collect a number of auxiliary results which are needed in order to prove the Neck Continuation Theorem and implement the surgery algorithm. The order has been arranged so as to make the consecutive choice of curvature thresholds and surgery parameters apparent. The proofs of these auxiliary results will be given in Sections \ref{proof.of.pseudoloc} -- \ref{proof.of.7.19}. 

We first establish a pseudolocality principle for the mean curvature flow. We begin with a definition.

\begin{definition}
\label{regular.flow}
Consider a ball $B$ in $\mathbb{R}^3$ and a one-parameter family of smooth surfaces $M_t \subset B$ such that $\partial M_t \subset \partial B$. Moreover, suppose that each surface $M_t$ bounds a domain $\Omega_t \subset B$. We say that the surfaces $M_t$ form a regular mean curvature flow if the surfaces $M_t$ form a smooth solution to mean curvature flow, except at finitely many times where one or more connected components of $\Omega_t$ may be removed.
\end{definition}

\begin{theorem}[Pseudolocality Principle] 
\label{pseudolocality}
There exist positive constants $\beta_0$ and $C$ such that the following holds. Suppose that $M_t$, $t \in [0,T]$, is a regular mean curvature flow in $B_4(0)$ in the sense of Definition \ref{regular.flow}. Moreover, we assume that the initial surface $M_0$ can be expressed as the graph of a (single-valued) function $u$ over a plane. If $\|u\|_{C^4} \leq \beta_0$, then 
\[|A|+|\nabla A|+|\nabla^2 A| \leq C\] 
for all $t \in [0,\beta_0] \cap [0,T]$ and all $x \in M_t \cap B_1(0)$.
\end{theorem}

Another important ingredient is the following curvature derivative estimate due to Haslhofer and Kleiner:

\begin{theorem}[cf. Haslhofer-Kleiner \cite{Haslhofer-Kleiner}, Theorem 1.8']
\label{interior.derivative.estimate}
Given any $\alpha \in (0,\frac{1}{1000}]$, there exists a constant $C(\alpha)$ with the following property. Suppose that $M_t$, $t \in [-1,0]$, is a regular mean curvature flow in the ball $B_4(0)$. Moreover, suppose that each surface $M_t$ is outward-minimizing within the ball $B_4(0)$. We further assume that the inscribed radius and the outer radius are at least $\frac{\alpha}{H}$ at each point on $M_t$. Finally, we assume that $M_0$ passes through the origin, and $H(0,0) \leq 1$. Then the surface $M_0$ satisfies $|\nabla A| \leq C(\alpha)$ and $|\nabla^2 A| \leq C(\alpha)$ at the origin.
\end{theorem}

In the following, we will fix an initial surface $M_0$ in $\mathbb{R}^3$. We assume that $M_0$ is closed, embedded, and has positive mean curvature. Moreover, let us fix a constant $\alpha \in (0,\frac{1}{1000}]$ such that the inscribed radius and the outer radius of the initial surface $M_0$ are at least $\frac{\alpha}{H}$. 

We next describe the necks on which we will perform surgery. 

\begin{definition}
Let $M$ be a mean convex surface in $\mathbb{R}^3$, and let $N$ be a region in $M$. As usual, we denote by $\nu$ the outward pointing unit normal vector field. We say that $N$ is an $(\hat{\alpha},\hat{\delta},\varepsilon,L)$-neck of size $r$ if (in a suitable coordinate system in $\mathbb{R}^3$) the following holds: 
\begin{itemize}
\item There exists a simple closed, convex curve $\Gamma \subset \mathbb{R}^2$ with the property that $\text{\rm dist}_{C^{20}}(r^{-1} \, N,\Gamma \times [-L,L]) \leq \varepsilon$.
\item At each point on $\Gamma$, the inscribed radius is at least $\frac{1}{(1+\hat{\delta}) \, \kappa}$, where $\kappa$ denotes the geodesic curvature of $\Gamma$.
\item We have $\sum_{l=1}^{18} |\nabla^l \kappa| \leq \frac{1}{100}$ at each point on $\Gamma$.
\item There exists a point on $\Gamma$ where the geodesic curvature $\kappa$ is equal to $1$.
\item The region $\{x + a \, \nu(x): x \in N, \, a \in (0,2\hat{\alpha} \, r)\}$ is disjoint from $M$.
\end{itemize}
\end{definition}

The last assumption is needed to ensure that, immediately after performing surgery, the resulting surface has outer radius at least $\frac{\alpha}{H}$ everywhere (cf. Proposition \ref{noncollapsing.preserved.under.surgery}). It turns out that the necks obtained via the Neck Detection Lemma satisfy this condition; see Theorem \ref{neck.detection.a} below.

Given an $(\hat{\alpha},\hat{\delta},\varepsilon,L)$-neck, we can perform surgery on $N$. The procedure depends on a parameter $\Lambda$, and will be explained in detail in Section \ref{construction.of.cap}. We will refer to this as $\Lambda$-surgery. The exact choice of $\Lambda$ will be specified later. 

\begin{theorem}[Properties of Surgery]
\label{properties.of.surgery}
Given any number $\hat{\alpha} > \alpha$, there exists a real number $\delta_0$ with the following significance. Suppose that we are given a pair of real numbers $\delta$ and $\hat{\delta}$ such that $\hat{\delta} < \delta < \delta_0$. Then we can find numbers $\bar{\varepsilon}$ and $\Lambda$, depending only on $\delta$ and $\hat{\delta}$, such that the following holds. Suppose that $N$ is an $(\hat{\alpha},\hat{\delta},\varepsilon,L)$-neck of size $r$ sitting in a mean convex surface in $\mathbb{R}^3$, where $\varepsilon \leq \bar{\varepsilon}$ and $\frac{L}{1000} \geq \Lambda$. If we perform a $\Lambda$-surgery on $N$, then the resulting surface $\tilde{N}$ will be $\frac{1}{1+\delta}$-noncollapsed. Furthermore, the outer radius is at least $\frac{\alpha}{H}$ at each point on $\tilde{N}$. Finally, if $\tilde{p} \in \tilde{N} \setminus N$ is a point in the surgically modified region, then either $\lambda_1(\tilde{p}) \geq 0$, or else there exists a point $p \in N$ such that $\lambda_1(\tilde{p}) \geq \lambda_1(p)$ and $H(\tilde{p}) \geq H(p)$.
\end{theorem}

A key point is that the deterioration in the noncollapsing constant can be made arbitrarily small by choosing $\varepsilon$ small and $\Lambda$ large.

\begin{assumption} 
\label{a.priori.assumptions}
In the following, we will assume that $M_t$ is a solution of the mean curvature flow which is interrupted by finitely many surgeries as in \cite{Huisken-Sinestrari3}, p.~145. We will assume that this flow satisfies the following assumptions: 
\begin{itemize}
\item The flow $M_t$ is smooth for $t \in [0,(100 \, \sup_{M_0} |A|)^{-2}]$.
\item Each surgery procedure involves performing a $\Lambda$-surgery on an $(\hat{\alpha},\hat{\delta},\varepsilon,L)$-neck of size $r \in [\frac{1}{2H_1},\frac{2}{H_1}]$, where $\hat{\alpha} > \alpha$, $\hat{\delta} \leq \frac{1}{10}$, $\frac{L}{1000} \geq \Lambda$, and $H_1 \geq \frac{(1000 \, \sup_{M_0} |A|)^2}{\inf_{M_0} H}$.
\item The region $\Omega_t$ enclosed by $M_t$ shrinks as $t$ increases.
\item For each $t$, the surface $M_t$ is outward-minimizing within the region $\Omega_0$.
\item For each $t$, the inscribed radius and the outer radius of $M_t$ are at least $\frac{\alpha}{H}$.
\end{itemize} 
\end{assumption}

The exact values of the parameters $\hat{\alpha}$, $\hat{\delta}$, $\Lambda$, $\varepsilon$, $L$, and $H_1$ will be specified later. 

In the first step, we want to apply the Pseudolocality Theorem to obtain derivative bounds shortly after a surgery. We begin by showing that surgeries are seperated in space:

\begin{proposition}[Separation of Surgery Regions]
\label{separation.of.surgery.regions}
Let $M_t$ be a mean curvature flow with surgery satisfying \ref{a.priori.assumptions}. Suppose that $t_0<t_1$ are two surgery times, and $x_0 \in M_{t_0+}$ and $x_1 \in M_{t_1+}$ are two points in the surgically modified regions. Then $|x_1-x_0| > \frac{\alpha}{1000} \, H_1^{-1}$. 
\end{proposition}

Thus, if $t_0$ is a surgery time and $x_0$ is a point in the surgically modified region at time $t_0+$, then the flow $M_t \cap B_{\frac{\alpha}{1000} \, H_1^{-1}}(x_0)$, $t > t_0$, is a regular flow in the sense of Definition \ref{regular.flow}. Using the Pseudolocality Theorem, we can draw the following conclusion:

\begin{proposition}
\label{consequence.of.pseudolocality}
There exist positive constants $\beta_* \in (0,\frac{\alpha}{1000})$ and $C_*$ with the following property. Let $M_t$ be a mean curvature flow with surgery satisfying Assumption \ref{a.priori.assumptions}. Suppose that $t_0$ is a surgery time and $x_0$ is a point in the surgically modified region at time $t_0+$. Then we have 
\[H_1^{-1} \, |A| + H_1^{-2} \, |\nabla A| + H_1^{-3} \, |\nabla^2 A| \leq C_*\] 
for all times $t \in (t_0,t_0+\beta_* \, H_1^{-2}]$ and all points $x \in M_t \cap B_{\beta_* \, H_1^{-1}}(x_0)$. The constants $\beta_*$ and $C_*$ may depend on the noncollapsing constant $\alpha$, but they do not depend on the surgery parameters $\hat{\alpha}$, $\hat{\delta}$, $\Lambda$, $\varepsilon$, $L$, and $H_1$.
\end{proposition}

The exact values of the surgery parameters will depend on the value of the constant in the derivative estimate, which in turn depends on $\beta_*$ and $C_*$. It is therefore critically important that the constants $\beta_*$ and $C_*$ do not depend on the exact choice of the surgery parameters $\hat{\alpha}$, $\hat{\delta}$, $\Lambda$, $\varepsilon$, $L$, and $H_1$.

Combining Proposition \ref{consequence.of.pseudolocality} with the interior gradient estimate of Haslhofer and Kleiner \cite{Haslhofer-Kleiner}, we obtain pointwise bounds for the first and second derivatives of the second fundamental form which hold even in the presence of surgeries.

\begin{proposition}[Pointwise Derivative Estimate]
\label{gradient.estimate}
There exists a constant $C_\#$ with the following significance. Suppose that $M_t$ is a mean curvature flow with surgery satisfying Assumption \ref{a.priori.assumptions}. Then $|\nabla A| \leq C_\# \, (H+H_1)^2$ and $|\nabla^2 A| \leq C_\# \, (H+H_1)^3$ for all times $t \geq (1000 \, \sup_{M_0} |A|)^{-2}$ and all points $x \in M_t$. The constant $C_\#$ may depend on the initial noncollapsing constant $\alpha$, but is independent of the surgery parameters $\hat{\alpha}$, $\hat{\delta}$, $\Lambda$, $\varepsilon$, $L$, and $H_1$.
\end{proposition}

Having fixed the constant $C_\#$ in the derivative estimate, we next define $\Theta = \frac{400}{\alpha}$, $\theta_0 = 10^{-6} \, \min\{\alpha,\frac{1}{C_\# \, \Theta^3}\}$, and $\hat{\alpha} = \frac{\alpha}{1-\frac{\theta_0}{8}}$. Hence, if we start at a point $(p_0,t_0)$ with $H(p_0,t_0) \geq \frac{H_1}{\Theta}$ and follow this point back in time, then the mean curvature at the resulting point will be between $\frac{1}{2} \, H(p_0,t_0)$ and $2 \, H(p_0,t_0)$, provided that $t \in (t_0-2\theta_0 \, H(p_0,t_0)^{-2},t_0]$. 

We next establish two auxiliary results concerning curves in the plane. It is here that we fix our choice of $\delta$ and $\hat{\delta}$. By applying these results to a blow-up limit that splits off line, we will show that the noncollapsing constants of a neck improve prior to surgery; this improvement offsets the deterioration of the noncollapsing constants under surgery (see Theorem \ref{properties.of.surgery} above).

\begin{proposition}
\label{choice.of.delta}
We can find a real number $\delta > 0$ such that the following holds: 
\begin{itemize} 
\item Suppose that $\Gamma$ is a (possibly non-closed) embedded curve in the plane with the property that $\kappa > 0$, $|\frac{d\kappa}{ds}| \leq C_\# \, (\kappa+2\Theta)^2$, and $|\frac{d^2\kappa}{ds^2}| \leq C_\# \, (\kappa+2\Theta)^3$. Moreover, suppose that the inscribed radius is at least $\frac{1}{(1+\delta) \, \kappa}$ at each point on $\Gamma$, and the outer radius is at least $\frac{\alpha}{\kappa}$ at each point on $\Gamma$. Finally, we assume that $\kappa(p)=1$ for some point $p \in \Gamma$. Then $L(\Gamma) \leq 3\pi$ and $\sup_\Gamma |\kappa-1| \leq \frac{1}{100}$. 
\item Suppose that $\Gamma_t$, $t \in (-2\theta_0,0]$, is a family of simple closed, convex curves in the plane which evolve by curve shortening flow. Assume that, for each $t \in (-2\theta_0,0]$, the curve $\Gamma_t$ satisfies the derivative estimates $|\frac{d\kappa}{ds}| \leq C_\# \, (\kappa+2\Theta)^2$ and $|\frac{d^2\kappa}{ds^2}| \leq C_\# \, (\kappa+2\Theta)^3$. Moreover, we assume that the inscribed radius is at least $\frac{1}{(1+\delta) \, \kappa}$ at each point on $\Gamma_t$, and the outer radius is at least $\frac{\alpha}{\kappa}$ at each point on $\Gamma_t$. Finally, we assume that the geodesic curvature of $\Gamma_0$ is equal to $1$ somewhere. Then the curve $\Gamma_0$ satisfies $\sum_{l=1}^{18} |\nabla^l \kappa| \leq \frac{1}{1000}$. Moreover, we have $\sup_{\Gamma_{-\theta_0}} \kappa \leq 1-\frac{\theta_0}{4}$.
\end{itemize}
\end{proposition}

We assume that $\delta$ is chosen sufficiently small so that $\delta < \delta_0$, where $\delta_0$ is the constant in Theorem \ref{properties.of.surgery}. In the next step, we choose $\hat{\delta}$ such that the following holds:

\begin{proposition}
\label{choice.of.hat.delta}
Given $\theta_0 > 0$ and $\delta > 0$, we can find a real number $\hat{\delta} \in (0,\delta)$ with the following property: Consider a simple closed, convex solution $\Gamma_t$, $t \in (-2\theta_0,0]$, of the curve shortening flow in the plane which satisfies the derivative estimates $|\frac{d\kappa}{ds}| \leq C_\# \, (\kappa+2\Theta)^2$ and $|\frac{d^2\kappa}{ds^2}| \leq C_\# \, (\kappa+2\Theta)^3$. Moreover, we assume that the inscribed radius is at least $\frac{1}{(1+\delta) \, \kappa}$ at each point on $\Gamma_t$, and the outer radius is at least $\frac{\alpha}{\kappa}$ at each point on $\Gamma_t$. Finally, we assume that the geodesic curvature of $\Gamma_0$ is equal to $1$ somewhere. Then $\Gamma_0$ is $\frac{1}{1+\hat{\delta}}$-noncollapsed.
\end{proposition}

Having fixed the values of $\alpha$, $\hat{\alpha}$, $\delta$, $\hat{\delta}$, we will choose $\bar{\varepsilon}$ and $\Lambda$ such that the conclusion of Theorem \ref{properties.of.surgery} holds. 

We next observe that the convexity estimates of Huisken and Sinestrari (cf. \cite{Huisken-Sinestrari1}, \cite{Huisken-Sinestrari2}) still hold for mean curvature flow with surgery.

\begin{proposition}[Huisken-Sinestrari \cite{Huisken-Sinestrari3}, Section 4]
\label{huisken.sinestrari.convexity}
Suppose that $\bar{\varepsilon}$ and $\Lambda$ are chosen such that the conclusion of Theorem \ref{properties.of.surgery} holds. Moreover, let $M_t$ be a mean curvature flow with surgery satisfying Assumption \ref{a.priori.assumptions}, where $\varepsilon \leq \bar{\varepsilon}$ and $L \geq 1000 \, \Lambda$. Given any $\eta > 0$, there exists a constant $C_1(\eta)$ such that $\lambda_1 \geq -\eta \, H - C_1(\eta)$. The constant $C_1(\eta)$ depends only on $\eta$ and the initial data, but is independent of the remaining surgery parameters $\varepsilon$, $L$, and $H_1$.
\end{proposition}

Theorem \ref{properties.of.surgery} implies that performing $\Lambda$-surgery on an $(\hat{\alpha},\hat{\delta},\varepsilon,L)$-neck will produce a surface which is $\frac{1}{1+\delta}$-noncollapsed, provided that $\varepsilon \leq \bar{\varepsilon}$ and $L \geq 1000 \, \Lambda$. This allows us to show that the cylindrical estimate from \cite{Brendle2} holds in the presence of surgeries:

\begin{proposition}[Cylindrical Estimate] 
\label{cylindrical.estimate}
Let $\delta$ and $\hat{\delta}$ be chosen as above. Moreover, suppose that $\bar{\varepsilon}$ and $\Lambda$ are chosen such that the conclusion of Theorem \ref{properties.of.surgery} holds. Finally, let $M_t$ be a mean curvature flow with surgery satisfying Assumption \ref{a.priori.assumptions}, where $\varepsilon \leq \bar{\varepsilon}$ and $L \geq 1000 \, \Lambda$. Then $\mu \leq (1+\delta) \, H + C \, H^{1-\sigma}$ for $0 \leq t \leq T$, where $\mu$ denotes the reciprocal of the inscribed radius. Here, $\sigma$ and $C$ may depend on $\delta$ and the initial data, but they are independent of the exact choice of $\varepsilon$, $L$, and $H_1$.
\end{proposition}

Using the convexity estimate and the cylindrical estimate, we are able to prove an analogue of the Neck Detection Lemma in \cite{Huisken-Sinestrari3}. In fact, we will need two different versions. 

\begin{theorem}[Neck Detection Lemma, Version A]
\label{neck.detection.a}
Let $\delta$ and $\hat{\delta}$ be chosen as above, and let $\bar{\varepsilon}$ and $\Lambda$ be chosen so that the conclusion of Theorem \ref{properties.of.surgery} holds. Let $M_t$ be a mean curvature flow with surgery satisfying Assumption \ref{a.priori.assumptions}, where $\varepsilon \leq \bar{\varepsilon}$ and $L \geq 1000 \, \Lambda$. Then, given $\varepsilon_0>0$ and $L_0 \geq 100$, we can find $\eta_0 > 0$ and $K_0$ with the following significance: Suppose that $t_0$ and $p_0 \in M_{t_0}$ satisfy 
\begin{itemize}
\item $H(p_0,t_0) \geq \max\{K_0,\frac{H_1}{\Theta}\}$, $\frac{\lambda_1(p_0,t_0)}{H(p_0,t_0)} \leq \eta_0$, 
\item the parabolic neighborhood $\hat{\mathcal{P}}(p_0,t_0,L_0+4,2\theta_0)$ does not contain surgeries.\footnote{See \cite{Huisken-Sinestrari3}, pp.~189--190, for the definition of $\hat{\mathcal{P}}(p_0,t_0,L_0+4,2\theta_0)$.}
\end{itemize}
Then $(p_0,t_0)$ lies at the center of an $(\hat{\alpha},\hat{\delta},\varepsilon_0,L_0)$-neck of size $H(p_0,t_0)^{-1}$. Finally, the constants $\eta_0$ and $K_0$ may depend on $\varepsilon_0$, $L_0$, $\delta$, $\hat{\delta}$, and the initial data, but they are independent of the remaining surgery parameters $\varepsilon$, $L$, and $H_1$.
\end{theorem}

\begin{theorem}[Neck Detection Lemma, Version B]
\label{neck.detection.b}
Let $\delta$ and $\hat{\delta}$ be chosen as above, and let $\bar{\varepsilon}$ and $\Lambda$ be chosen so that the conclusion of Theorem \ref{properties.of.surgery} holds. Let $M_t$ be a mean curvature flow with surgery satisfying Assumption \ref{a.priori.assumptions}, where $\varepsilon \leq \bar{\varepsilon}$ and $L \geq 1000 \, \Lambda$. Then, given $\theta$, $\varepsilon_0>0$ and $L_0 \geq 100$, we can find positive numbers $\eta_0$ and $K_0$ with the following significance: Suppose that $t_0$ and $p_0 \in M_{t_0}$ satisfy 
\begin{itemize}
\item $H(p_0,t_0) \geq \max\{K_0,\frac{H_1}{\Theta}\}$, $\frac{\lambda_1(p_0,t_0)}{H(p_0,t_0)} \leq \eta_0$, 
\item the parabolic neighborhood $\hat{\mathcal{P}}(p_0,t_0,L_0+4,\theta)$ does not contain surgeries.
\end{itemize}
Let us dilate the surface $\{x \in M_{t_0}: d_{g(t_0)}(p_0,x) \leq L_0 \, H(p_0,t_0)^{-1}\}$ by the factor $H(p_0,t_0)$. Then the resulting surface is $\varepsilon_0$-close to a product $\Gamma \times [-L_0,L_0]$ in the $C^3$-norm. Here, $\Gamma$ is a closed, convex curve satisfying $L(\Gamma) \leq 3\pi$ and $\sup_\Gamma |\kappa-1| \leq \frac{1}{100}$. The constant $K_0$ may depend on $\theta$, $\varepsilon_0$, $L_0$, $\delta$, $\hat{\delta}$, and the initial data, but they are independent of the remaining surgery parameters $\varepsilon$, $L$, and $H_1$.
\end{theorem}

The proof of the Neck Continuation Theorem in Section \ref{neck.contin} will require both versions of the Neck Detection Lemma. We will describe the proof of Version A in Section \ref{proof.of.neck.detection.lemma.version.a}. (The proof of Version B is analogous.) The main difference between the two versions is that Version A requires the assumption that $\hat{\mathcal{P}}(p_0,t_0,L_0+4,2\theta_0)$ does not contain surgeries, whereas Version B only requires that the parabolic neighborhood $\hat{\mathcal{P}}(p_0,t_0,L_0+4,\theta)$ is free of surgeries. (Note that $\theta$ can be much smaller than $\theta_0$.)

The following next result serves as a replacement for Lemma 7.12 in \cite{Huisken-Sinestrari3}: 

\begin{proposition}[Replacement for Lemma 7.12 in \cite{Huisken-Sinestrari3}]
\label{7.12}
Let $M_t$ be a mean curvature flow with surgery satisfying Assumption \ref{a.priori.assumptions}. Suppose that $(p_1,t_1)$ is a point in spacetime such that $H(p_1,t_1) \geq H_1$ and the parabolic neighborhood $\hat{\mathcal{P}}(p_1,t_1,\tilde{L}+4,2\theta_0)$ contains at least one point belonging to a surgery region. Then there exists a point $q_1 \in M_{t_1}$ and an open set $V \subset M_{t_1}$ such that 
$d_{g(t_1)}(p_1,q_1) \leq (\tilde{L}+4) \, H(p_1,t_1)^{-1}$, $\{x \in M_{t_1}: d_{g(t_1)}(q_1,x) \leq 500 \, H_1^{-1}\} \subset V$, and $V$ is diffeomorphic to a disk. Moreover, the mean curvature is at most $40 \, H_1$ at each point in $V$.
\end{proposition}

To construct $V$, we consider a surgical cap that was inserted shortly before time $t_1$. We then follow this cap forward in time (see Section \ref{proof.of.7.12} below).

Since we have a bound for the gradient of the mean curvature, we can apply Theorem 7.14 in \cite{Huisken-Sinestrari3}. This gives the following result: 

\begin{proposition}[Huisken-Sinestrari \cite{Huisken-Sinestrari3}, Theorem 7.14]
\label{7.14}
Consider a closed surface in $\mathbb{R}^3$ which satisfies the estimate $|\nabla A| \leq C_\# \, (H+H_1)^2$ for suitable constants $C_\#$ and $H_1$. Then, given any $\eta>0$, we can find large numbers $\rho$ and $\gamma_0$ (depending only on $C_\#$ and $\eta$) with the following significance. Suppose that $p$ is a point on the surface with $\lambda_1(p) > \eta \, H(p)$ and $H(p) \geq \gamma_0 \, H_1$. Then either $\lambda_1 > \eta \, H > 0$ everywhere on the surface, or else there exists a point $q$ such that $\lambda_1(q) \leq \eta \, H(q)$; $d(p,q) \leq \frac{\rho}{H(p)}$; and $H(q') \geq \frac{H(p)}{\gamma_0} \geq H_1$ for all points $q'$ satisfying $d(p,q') \leq \frac{\rho}{H(p)}$. In particular, $H(q) \geq \frac{H(p)}{\gamma_0} \geq H_1$.
\end{proposition}

Moreover, using the noncollapsing property we can prove an analogue of Lemma 7.19 in \cite{Huisken-Sinestrari3}. This result will be needed for the proof of the Neck Continuation Theorem. 

\begin{proposition}[Replacement for Lemma 7.19 in \cite{Huisken-Sinestrari3}]
\label{7.19}
Let $\Sigma$ be an embedded surface in $\mathbb{R}^3$ which is $\alpha$-noncollapsed, and let $y_1<y_2$ be two real numbers. We assume that the surface $\Sigma$ is contained in the cylinder $\{(x_1,x_2,x_3) \in \mathbb{R}^3: x_1^2+x_2^2 \leq 100, \, y_1 \leq x_3 \leq y_2\}$. Moreover, we assume that $\partial \Sigma = \Gamma_1 \cup \Gamma_2$, where $\Gamma_1 \subset \{x \in \mathbb{R}^3: x_3=y_1\}$ and $\Gamma_2 \subset \{x \in \mathbb{R}^3: x_3=y_2\}$. Then we have $H(x) \geq \frac{4}{\Theta}$ for all points $x \in \Sigma$ satisfying $x_3 \in [y_1+1,y_2]$ and $\langle \nu(x),e_3 \rangle \geq 0$. Here, $\nu$ denotes the outward-pointing unit normal to $\Sigma$ and $\Theta = \frac{400}{\alpha}$.
\end{proposition}

\section{The Neck Continuation Theorem and the proof of Theorem \ref{main.theorem}}

\label{neck.contin}

In this section, we use the auxiliary results collected in Section \ref{overview} to establish an analogue of the Neck Continuation Theorem of Huisken and Sinestrari \cite{Huisken-Sinestrari3}. 

We begin by finalizing our choice of the surgery parameters. This step is similar to the discussion on pp.~208--209 in \cite{Huisken-Sinestrari3}. Recall that the parameters $\delta$, $\hat{\delta}$, $\hat{\alpha}$ and the constants $C_\#$, $\theta_0$, $\Theta$ have already been chosen at this stage. Moreover, we have chosen $\bar{\varepsilon}$ and $\Lambda$ so that the conclusion of Theorem \ref{properties.of.surgery} holds. 

In the next step, we choose numbers $\varepsilon_0$ and $L_0$ so that $\varepsilon_0 < \bar{\varepsilon}$ and $L_0 > 1000 \, \Lambda$. In addition, we require that the mean curvature on an $(\hat{\alpha},\hat{\delta},\varepsilon_0,L_0)$-neck varies by at most a factor of $1+L_0^{-1}$. (This can always be achieved by choosing $\varepsilon_0$ very small.) We then choose real numbers $\eta_0 > 0$ and $K_0 > 1000 \, \sup_{M_0} |A|$ so that the conclusion of Version A of the Neck Detection Lemma can be applied for each $\tilde{L} \in [100,L_0]$. In other words, if $(p_0,t_0)$ satisfies $H(p_0,t_0) \geq \max \{K_0,\frac{H_1}{\Theta}\}$, $\lambda_1(p_0,t_0) \leq \eta_0 \, H(p_0,t_0)$, and if the parabolic neighborhood $\hat{\mathcal{P}}(p_0,t_0,\tilde{L}+4,2\theta_0)$ is free of surgeries for some $\tilde{L} \in [100,L_0]$, then $(p_0,t_0)$ lies at the center of a $(\hat{\alpha},\hat{\delta},\varepsilon_0,\tilde{L})$-neck in $M_{t_0}$.

In the next step, we put $\varepsilon_1 = \frac{\eta_0}{10}$. By Version A of the Neck Detection Lemma, we can find constants $\eta_1<\eta_0$ and $K_1>K_0$ such that the following holds: if $(p_0,t_0)$ satisfies $H(p_0,t_0) \geq \max \{K_1,\frac{H_1}{\Theta}\}$, $\lambda_1(p_0,t_0) \leq \eta_1 \, H(p_0,t_0)$, and if the parabolic neighborhood $\hat{\mathcal{P}}(p_0,t_0,104,2\theta_0)$ is free of surgeries, then $(p_0,t_0)$ lies at the center of a $(\hat{\alpha},\hat{\delta},\varepsilon_1,100)$-neck in $M_{t_0}$. 

Having chosen $\eta_1$, we next choose $\gamma_0$ and $\rho$ so that the conclusion of Proposition \ref{7.14} holds with $\eta=\eta_1$.

By Version B of the Neck Detection Lemma, we can find a number $K_2>K_1$ such that the following holds: Suppose that $(p_0,t_0)$ satisfies $H(p_0,t_0) \geq \max\{K_2,\frac{H_1}{\Theta}\}$ and  $\lambda_1(p_0,t_0) \leq 0$, and that the parabolic neighborhood $\hat{\mathcal{P}}(p_0,t_0,104,10^{-6} \, \Theta^{-2} \, \gamma_0^{-2})$ does not contain surgeries. Then, if we dilate the surface $\{x \in M_{t_0}: d_{g(t_0)}(p_0,x) \leq 100 \, H(p_0,t_0)^{-1}\}$ by the factor $H(p_0,t_0)$, the resulting surface is $\frac{\varepsilon_1}{10}$-close to a product $\Gamma \times [-100,100]$ in the $C^3$-norm. Here, $\Gamma$ is a closed, convex curve satisfying $L(\Gamma) \leq 3\pi$ and $\sup_\Gamma |\kappa-1| \leq \frac{1}{100}$. 

Finally, we choose $H_1 \geq 1000 \, \Theta \, K_2$, and define $H_2 = 1000 \, \gamma_0 \, H_1$, and $H_3 = 10 \, H_2$. %We note that the distance between consecutive surgery times is bounded from below by 
%\[\frac{1}{3} \, (H_2^{-2} - H_3^{-2}) > \frac{1}{4} \, (1000 \, \gamma_0 \, H_1)^{-2} = 10^{-6} \, \Theta^{-2} \, \gamma_0^{-2} \, \Big ( \frac{2H_1}{\Theta} \Big )^{-2}\] 
%(compare \cite{Huisken-Sinestrari3}, p.~210). 

Recall that the Neck Detection Lemma requires that a certain parabolic neighborhood is free of surgeries. It turns out that this assumption is not needed when the curvature is at least $1000 \, H_1$:

\begin{proposition}
\label{7.10}
Suppose that $M_t$ is a mean curvature flow with surgeries satisfying Assumption \ref{a.priori.assumptions}, where $\varepsilon \leq \bar{\varepsilon}$ and $L \geq 1000 \, \Lambda$. Moreover, suppose that $(p_0,t_0)$ satisfies $H(p_0,t_0) \geq 1000 \, H_1$ and $\lambda_1(p_0,t_0) \leq \eta_0 \, H(p_0,t_0)$, where $\eta_0$ and $H_1$ are defined as above. Then $p_0$ lies at the center of an $(\hat{\alpha},\hat{\delta},\varepsilon_0,L_0)$-neck.
\end{proposition}

\textbf{Proof.} 
We distinguish two cases: 

\textit{Case 1:} Suppose first that the parabolic neighborhood $\hat{\mathcal{P}}(p_0,t_0,104,2\theta_0)$ contains a point modified by surgery. By Proposition \ref{7.12}, we can find a point $q \in M_{t_0}$ and an open set $V \subset \{x \in M_{t_0}: H(x,t_0) \leq 40 \, H_1\}$ such that $d_{g(t_0)}(p_0,q) \leq 104 \, H(p_0,t_0)^{-1}$ and $\{x \in M_{t_0}: d_{g(t_0)}(q,x) \leq 500 \, H_1^{-1}\} \subset V$. Clearly, $p_0 \in V$. Consequently, $H(p_0,t_0) \leq 40 \, H_1$, contrary to our assumption. 

\textit{Case 2:} We now assume that the parabolic neighborhood $\hat{\mathcal{P}}(p_0,t_0,104,2\theta_0)$ is free of surgeries. Let $\tilde{L} \in [100,L_0]$ be the largest number with the property that $\hat{\mathcal{P}}(p_0,t_0,\tilde{L}+4,2\theta_0)$ is free of surgeries. By Version A of the Neck Detection Lemma, the point $(p_0,t_0)$ lies at the center of an $(\hat{\alpha},\hat{\delta},\varepsilon_0,\tilde{L})$-neck $N$. If $\tilde{L}=L_0$, we are done. Hence, it remains to consider the case when $\tilde{L}<L_0$. In this case, the parabolic neighborhood $\hat{\mathcal{P}}(p_0,t_0,\tilde{L}+5,2\theta_0)$ must contain a point modified by surgery. By Proposition \ref{7.12}, we can find a point $q \in M_{t_0}$ and an open set $V \subset \{x \in M_{t_0}: H(x,t_0) \leq 40 \, H_1\}$ such that $d_{g(t_0)}(p_0,q) \leq (\tilde{L}+5) \, H(p_0,t_0)^{-1}$ and $\{x \in M_{t_0}: d_{g(t_0)}(q,x) \leq 500 \, H_1^{-1}\} \subset V$. Since the set $\{x \in M_{t_0}: d_{g(t_0)}(p_0,x) \leq (\tilde{L}-1) \, H(p_0,t_0)^{-1}\}$ is contained in $N$, we conclude that $\text{\rm dist}_{g(t_0)}(q,N) \leq 6 \, H(p_0,t_0)^{-1} \leq 6 \, H_1^{-1}$. Consequently, we have $N \cap V \neq \emptyset$. On the other hand, we have $H \geq \frac{1}{2} \, H(p_0,t_0) \geq 500 \, H_1$ at each point on $N$ and $H \leq 40 \, H_1$ at each point on $V$. This is a contradiction. This completes the proof of Proposition \ref{7.10}. \\

The following result is the analogue of the Neck Continuation Theorem in \cite{Huisken-Sinestrari3}:

\begin{theorem}[Neck Continuation Theorem]
\label{neck.continuation}
Suppose that $M_t$ is a mean curvature flow with surgery satisfying Assumption \ref{a.priori.assumptions}, where $\varepsilon \leq \bar{\varepsilon}$ and $L \geq 1000 \, \Lambda$. Suppose that $(p_0,t_0)$ satisfies $H(p_0,t_0) \geq 1000 \, H_1$ and $\lambda_1(p_0,t_0) \leq \eta_1 \, H(p_0,t_0)$, where $\eta_1$ and $H_1$ are defined as above. Then there exists a finite collection of points $p_1,\hdots,p_l$ with the following properties: 
\begin{itemize}
\item For each $i=0,1,\hdots,l$, the point $p_i$ lies at the center of an $(\hat{\alpha},\hat{\delta},\varepsilon_0,L_0)$-neck $N^{(i)} \subset M_{t_0}$, and we have $H(p_i,t_0) \geq H_1$.
\item For each $i=1,\hdots,l-1$, the point $p_{i+1}$ lies on the neck $N^{(i)}$, and we have $\text{\rm dist}_{g(t_0)}(p_{i+1},\partial N^{(i)} \setminus N^{(i-1)}) \in [(L_0-100) \, H(p_i,t_0)^{-1},(L_0-50) \, H(p_i,t_0)^{-1}]$. 
\item Finally, at least one of the following four statements holds: either the union $\mathcal{N} = \bigcup_{i=1}^l N^{(i)}$ covers the entire surface; or $H(p_l,t_0) \in [H_1,2H_1]$; or there exists a closed curve in $\mathcal{N} \cap \{x \in M_{t_0}: H(x,t_0) \leq 40 \, H_1\}$ which is homotopically non-trivial in $\mathcal{N}$ and bounds a disk in $\{x \in M_{t_0}: H(x,t_0) \leq 40 \, H_1\}$; or the outer boundary $\partial N^{(k)} \setminus N^{(k-1)}$ bounds a convex cap.
\end{itemize}
\end{theorem}

We now describe the proof of the Neck Continuation Theorem. Most of the arguments in \cite{Huisken-Sinestrari3} carry over to our situation. However, the proof of Lemma 7.19 does not work in our setting. The reason is that the gradient estimate in \cite{Huisken-Sinestrari3} works on all scales, whereas the gradient estimate in Proposition \ref{gradient.estimate} becomes weaker when the curvature is much smaller than $H_1$. We will use Proposition \ref{7.19} to overcome this problem. 

\textbf{Proof.} By Proposition \ref{7.10}, the point $p_0$ lies at the center of an $(\hat{\alpha},\hat{\delta},\varepsilon_0,L_0)$-neck $N^{(0)} \subset M_{t_0}$. The construction of the points $p_1,p_2,\hdots$ is by induction. Suppose that we have constructed points $p_1,\hdots,p_k$ and necks $N^{(1)},\hdots,N^{(k)}$ with the following properties: 
\begin{itemize}
\item For each $i=0,1,\hdots,k$, the point $p_i$ lies at the center of an $(\hat{\alpha},\hat{\delta},\varepsilon_0,L_0)$-neck $N^{(i)} \subset M_{t_0}$, and we have $H(p_i,t_0) \geq H_1$.
\item For each $i=1,\hdots,k-1$, the point $p_{i+1}$ lies on the neck $N^{(i)}$, and we have $\text{\rm dist}_{g(t_0)}(p_{i+1},\partial N^{(i)} \setminus N^{(i-1)}) \in [(L_0-100) \, H(p_i,t_0)^{-1},(L_0-50) \, H(p_i,t_0)^{-1}]$. 
\end{itemize}
If $H(p_k,t_0) \in [H_1,2H_1]$, then we are done. Hence, for the remainder of the proof, we will assume that $H(p_k,t_0) \geq 2H_1$. We break the discussion into several cases:

\textit{Case 1:} Suppose that the there exists a point $p \in N^{(k)}$ such that $\text{\rm dist}_{g(t_0)}(p,\partial N^{(k)} \setminus N^{(k-1)}) \in [(L_0-100) \, H(p_k,t_0)^{-1},(L_0-50) \, H(p_k,t_0)^{-1}]$, and the parabolic neighborhood $\hat{\mathcal{P}}(p,t_0,L_0+4,2\theta_0)$ contains a point modified by surgery. In this case, Proposition \ref{7.12} implies that there exists a point $q \in M_{t_0}$ and an open set $V \subset \{x \in M_{t_0}: H(x,t_0) \leq 40 \, H_1\}$ such that $d_{g(t_0)}(p,q) \leq (L_0+4) \, H(p,t_0)^{-1}$, $\{x \in M_{t_0}: d_{g(t_0)}(q,x) \leq 500 \, H_1^{-1}\} \subset V$, and $V$ is diffeomorphic to a disk. 

By our choice of $\varepsilon_0$ and $L_0$, the mean curvature on $N^{(k)}$ varies at most by a factor $1+L_0^{-1}$. Hence, $H(p_k,t_0) \leq (1+L_0^{-1}) \, H(p,t_0)$. Since the set $\{x \in M_{t_0}: d_{g(t_0)}(p,x) \leq (L_0-100) \, H(p_k,t_0)^{-1}\}$ is contained in $N^{(k)}$, we conclude that 
\begin{align*} 
\text{\rm dist}_{g(t_0)}(q,N^{(k)}) 
&\leq (L_0+4) \, H(p,t_0)^{-1}-(L_0-100) \, H(p_k,t_0)^{-1} \\ 
&\leq (L_0+4) \, (1+L_0^{-1}) \, H(p_k,t_0)^{-1}-(L_0-100) \, H(p_k,t_0)^{-1} \\ 
&\leq 200 \, H(p_k,t_0)^{-1} \\ 
&\leq 100 \, H_1^{-1}. 
\end{align*} 
Consequently, there exists a closed curve which is contained in $N^{(k)} \cap V$ and is homotopically non-trivial in $N^{(k)}$. Since $V$ is diffeomorphic to a disk, this curve bounds a disk in $V$, and we are done.

\textit{Case 2:} We now assume that the parabolic neighborhood $\hat{\mathcal{P}}(p,t_0,L_0+4,2\theta_0)$ is free of surgeries for all points $p \in N^{(k)}$ satisfying $\text{\rm dist}_{g(t_0)}(p,\partial N^{(k)} \setminus N^{(k-1)}) \in [(L_0-100) \, H(p_k,t_0)^{-1},(L_0-50) \, H(p_k,t_0)^{-1}]$. There are two possibilities now: 

\textit{Subcase 2.1:} Suppose that there exists a point $p \in N^{(k)}$ with the property that $\text{\rm dist}_{g(t_0)}(p,\partial N^{(k)} \setminus N^{(k-1)}) \in [(L_0-100) \, H(p_k,t_0)^{-1},(L_0-50) \, H(p_k,t_0)^{-1}]$ and $\lambda_1(p,t_0) \leq \eta_0 \, H(p,t_0)$. By Version A of the Neck Detection Lemma, the point $p$ lies at the center of an $(\hat{\alpha},\hat{\delta},\varepsilon_0,L_0)$-neck $N$. Moreover, since $p \in N^{(k)}$ and $H(p_k,t_0) \geq 2H_1$, we have $H(p,t_0) \geq H_1$. Hence, we can put $p^{(k+1)} := p$ and $N^{(k+1)} := N$ and continue the process. 

\textit{Subcase 2.2:} Suppose that $\lambda_1(p,t_0) > \eta_0 \, H(p,t_0)$ for all points $p \in N^{(k)}$ satisfying $\text{\rm dist}_{g(t_0)}(p,\partial N^{(k)} \setminus N^{(k-1)}) \in [(L_0-100) \, H(p_k,t_0)^{-1},(L_0-50) \, H(p_k,t_0)^{-1}]$. Let $\mathcal{N} = \bigcup_{i=0}^k N^{(i)}$, and let $\mathcal{A}$ be the set of all points $x \in \mathcal{N}$ satisfying $\text{\rm dist}_{g(t_0)}(p,\partial N^{(k)} \setminus N^{(k-1)}) \geq (L_0-50) \, H(p_k,t_0)^{-1}$ and $\lambda_1(x,t_0) \leq \eta_1 \, H(x,t_0)$. The assumptions of Theorem \ref{neck.continuation} imply that the initial point $p_0$ belongs to $\mathcal{A}$, so $\mathcal{A}$ is non-empty. Let us consider a point $p^*$ which has maximal intrinsic distance from $p_0$ among all points in $\mathcal{A}$.

\textit{Subcase 2.2.1:} Suppose that the parabolic neighborhood $\hat{\mathcal{P}}(p^*,t_0,104,2\theta_0)$ contains a point modified by surgery. In this case, Proposition \ref{7.12} implies that there exists a point $q \in M_{t_0}$ and an open set $V \subset \{x \in M_{t_0}: H(x,t_0) \leq 40 \, H_1\}$ such that $d_{g(t_0)}(p^*,q) \leq 104 \, H(p^*,t_0)^{-1}$, $\{x \in M_{t_0}: d_{g(t_0)}(q,x) \leq 500 \, H_1^{-1}\} \subset V$, and $V$ is diffeomorphic to a disk. Since $H(p^*,t_0) \geq \frac{H_1}{2}$, this implies 
\begin{align*} 
&\{x \in M_{t_0}: d_{g(t_0)}(p^*,x) \leq 100 \, H(p^*,t_0)^{-1}\} \\ 
&\subset \{x \in M_{t_0}: d_{g(t_0)}(q,x) \leq 204 \, H(p^*,t_0)^{-1}\} \\ 
&\subset \{x \in M_{t_0}: d_{g(t_0)}(q,x) \leq 500 \, H_1^{-1}\} \\ 
&\subset V. 
\end{align*}
Consequently, there exists a closed curve in $\mathcal{N} \cap V$ which is homotopically non-trivial in $\mathcal{N}$. This curve bounds a disk which is contained in $V$. Hence, we can again terminate the process.

\textit{Subcase 2.2.2:} Suppose, finally, that the parabolic neighborhood $\hat{\mathcal{P}}(p^*,t_0,104,2\theta_0)$ is free of surgeries. In this case, Version A of the Neck Detection Lemma implies that the point $p^*$ lies at the center of an $(\hat{\alpha},\hat{\delta},\varepsilon_1,100)$-neck $N^*$. Clearly, $\lambda_1 \leq \varepsilon_1 \, H$ at each point on $N^*$. Consequently, the set $N^*$ is disjoint from the set $\{p \in N^{(k)}: \text{\rm dist}_{g(t_0)}(p,\partial N^{(k)} \setminus N^{(k-1)}) \in [(L_0-100) \, H(p_k,t_0)^{-1},(L_0-50) \, H(p_k,t_0)^{-1}]\}$. Furthermore, since $p^*$ has maximal distance from $p_0$ among all points in $\mathcal{A}$, we conclude that the part of $\mathcal{N}$ that lies between the neck $N^*$ and the set $\{p \in N^{(k)}: \text{\rm dist}_{g(t_0)}(p,\partial N^{(k)} \setminus N^{(k-1)}) \in [(L_0-100) \, H(p_k,t_0)^{-1},(L_0-50) \, H(p_k,t_0)^{-1}]\}$ is strictly convex. 

Let $\omega$ be a unit vector in $\mathbb{R}^3$ which is parallel to the axis of the neck $N^*$. The arguments on p.~214 of \cite{Huisken-Sinestrari3} imply that $\langle \nu,\omega \rangle \geq -\varepsilon_1$ for all points $p \in N^{(k)}$ satisfying $\text{\rm dist}_{g(t_0)}(p,\partial N^{(k)} \setminus N^{(k-1)}) \in [(L_0-100) \, H(p_k,t_0)^{-1},(L_0-50) \, H(p_k,t_0)^{-1}]$. Moreover, we have $\lambda_1(p,t_0) > \eta_0 \, H(p,t_0)$ for all points $p \in N^{(k)}$ satisfying $\text{\rm dist}_{g(t_0)}(p,\partial N^{(k)} \setminus N^{(k-1)}) \in [(L_0-100) \, H(p_k,t_0)^{-1},(L_0-50) \, H(p_k,t_0)^{-1}]$. Putting these facts together (and using the fact that $\eta_0 \geq 10\varepsilon_1$), we conclude that $\langle \nu,\omega \rangle \geq 4\varepsilon_1$ for all points $p \in N^{(k)}$ satisfying $\text{\rm dist}_{g(t_0)}(p,\partial N^{(k)} \setminus N^{(k-1)}) \in [(L_0-100) \, H(p_k,t_0)^{-1},(L_0-75) \, H(p_k,t_0)^{-1}]$. 

We claim that the boundary curve $\partial N^{(k)} \setminus N^{(k-1)}$ bounds a convex cap. To prove this, we follow the argument on pp.~215-216 of \cite{Huisken-Sinestrari3}. Let us choose a curve $\Gamma_0$ such that $\Gamma_0 \subset \{p \in N^{(k)}: \text{\rm dist}_{g(t_0)}(p,\partial N^{(k)} \setminus N^{(k-1)}) \in [(L_0-100) \, H(p_k,t_0)^{-1},(L_0-75) \, H(p_k,t_0)^{-1}]\}$ and $\Gamma_0$ is contained in a plane orthogonal to $\omega$. For each point on $\Gamma_0$, we solve the ODE $\dot{\gamma} = \frac{\omega^T(\gamma)}{|\omega^T(\gamma)|^2}$, where $\omega^T(\gamma)$ denotes the projection of $\omega$ to the tangent plane to $M_{t_0}$ at the point $\gamma$. This gives a family of curves $\Gamma_y \subset M_{t_0}$, each of which is contained in a plane orthogonal to $\omega$. The curves $\Gamma_y$ are well-defined for $y \in [0,y_{\text{\rm max}})$. Moreover, there exists a point $p \in \Gamma_0$ such that $\nu(\gamma(y,p)) \to \omega$ as $y \to y_{\text{\rm max}}$. 

Following the arguments on p.~215 in \cite{Huisken-Sinestrari3}, we can show that the inequalities 
\begin{equation} 
\tag{$\star$}
\langle \nu,\omega \rangle < 1, \qquad \lambda_1 > 0, \qquad H > \frac{2H_1}{\Theta}, \qquad \langle \nu,\omega \rangle > \varepsilon_1 
\end{equation}
hold for all $y \in [0,y_{\text{\rm max}})$. Indeed, the inequalities in ($\star$) are clearly satisfied for $y=0$. If one of the inequalities in ($\star$) fails for some $y>0$, we consider the smallest value of $y$ for which that happens. The first inequality in ($\star$) cannot fail first by definition of $y_{\text{\rm max}}$. If the second inequality in ($\star$) is the first one to fail, then we have $\lambda_1=0$. Since $H \geq \frac{2H_1}{\Theta}$, we may apply Version B of the Neck Detection Lemma to conclude that we are $\frac{\varepsilon_1}{10}$-close to a Cartesian product, but this is ruled out by the fourth inequality in ($\star$). If the third inequality in ($\star$) is the first one to fail, we obtain a contradiction with Proposition \ref{7.19}. Finally, $\langle \nu,\omega \rangle$ is montone increasing in $y$ as long as $\lambda_1$ remains nonnegative; this implies that the fourth inequality in ($\star$) cannot fail first. Thus, the inequalities in ($\star$) hold for all $y \in [0,y_{\text{\rm max}})$. Consequently, the union of the curves $\Gamma_y$ is a convex cap, and we can terminate the process. This completes the construction of the sequence $p_1,p_2,\hdots$.

If the sequence $p_1,p_2,\hdots$ terminates after finitely many steps, then the theorem is proved. On the other hand, if the sequence $p_1,p_2,\hdots$ never terminates, then the necks $N^{(1)},N^{(2)},\hdots$ will eventually cover the entire surface. This completes the proof of Theorem \ref{neck.continuation}. \\

Having established the Neck Continuation Theorem, we can now implement the surgery algorithm of Huisken and Sinestrari \cite{Huisken-Sinestrari3}, and complete the proof of Theorem \ref{main.theorem}. Starting from the given initial surface $M_0$, we run the mean curvature flow until the maximum of the mean curvature reaches the threshold $H_3$ for the first time. Let us denote this time by $T_1$. By a result of Andrews \cite{Andrews}, the inscribed radius and the outer radius are bounded from below by $\frac{\alpha}{H}$ for $0 \leq t \leq T_1$. Moreover, it is easy to see that the surfaces $M_t$ are outward-minimizing for $0 \leq t \leq T_1$. Therefore, Assumption \ref{a.priori.assumptions} is satisfied for $0 \leq t \leq T_1$. Consequently, we may apply the Neck Detection Lemma and the Neck Continuation Theorem for $0 \leq t \leq T_1$. By performing surgery on suitably chosen $(\hat{\alpha},\hat{\delta},\varepsilon_0,L_0)$-necks at time $T_1$, we can remove all regions where the mean curvature is between $H_2$ and $H_3$. Hence, immediately after surgery, the maximum of the mean curvature drops to a level below $H_2$. We then run the flow again until the maximum of the mean curvature reaches $H_3$ for the second time. Let us denote this time by $T_2$. We claim that, for $0 \leq t \leq T_2$, the flow satisfies Assumption \ref{a.priori.assumptions} with $\varepsilon=\varepsilon_0$ and $L=L_0$. Indeed, Theorem \ref{properties.of.surgery} implies that the inscribed radius and the outer radius of the surface $M_{T_1+}$ are bounded from below by $\frac{\alpha}{H}$, and this property continues to hold for all $T_1 < t \leq T_2$ by a result of Andrews \cite{Andrews}. Furthermore, the outward-minimizing property follows from work of Head (see \cite{Head}, Lemma 5.2). Therefore, Assumption \ref{a.priori.assumptions} is satisfied for $0 \leq t \leq T_2$ with $\varepsilon=\varepsilon_0$ and $L=L_0$. Hence, we can apply the Neck Detection Lemma and the Neck Continuation Theorem for $0 \leq t \leq T_2$. By performing surgery on suitably chosen $(\hat{\alpha},\hat{\delta},\varepsilon,L)$-necks, we can push the maxmimum of the mean curvature below $H_2$. We then restart the flow again. This process can be repeated until the solution becomes extinct.

\section{Proof of the pseudolocality principle (Theorem \ref{pseudolocality})}

\label{proof.of.pseudoloc} 

We first recall the following analogue of Shi's local derivative estimate for the Ricci flow. The argument given here is standard and follows the proof in Ecker-Huisken \cite{Ecker-Huisken2}; see also \cite{Ecker}, Proposition 3.22.

\begin{lemma}
\label{shi.1}
Suppose that $M_t$, $t \in [0,T]$, is a regular mean curvature flow in $B_4(0)$ in the sense of Definition \ref{regular.flow}. Moreover, we assume that $|A| \leq 1$ for all $t \in [0,T]$ and all $x \in M_t \cap B_4(0)$. Finally, we assume that $|\nabla A| \leq 1$ for all $x \in M_0 \cap B_4(0)$. Then $|\nabla A| \leq C$ for all $t \in [0,1] \cap [0,T]$ and all $x \in M_t \cap B_2(0)$.
\end{lemma}

\textbf{Proof.} 
Consider the cutoff function $\psi(x) = 1-\frac{|x|^2}{16}$. A straightforward calculation gives 
\[\frac{\partial}{\partial t} (\psi^2 \, |\nabla A|^2) \leq \Delta (\psi^2 \, |\nabla A|^2) + C_0 \, |\nabla A|^2\] 
for all $t \in [0,T]$ and all $x \in M_t \cap B_4(0)$. This implies that 
\[\frac{\partial}{\partial t} (\psi^2 \, |\nabla A|^2 + C_0 \, |A|^2) \leq \Delta (\psi^2 \, |\nabla A|^2 + C_0 \, |A|^2) + C_1\] 
for all $t \in [0,T]$ and all $x \in M_t \cap B_4(0)$. Applying the maximum principle to the function $\psi^2 \, |\nabla A|^2 + C_0 \, |A|^2 - C_1 \, t$, we obtain 
\begin{align*} 
&\sup_{t \in [0,T]} \sup_{x \in M_t \cap B_4(0)} (\psi^2 \, |\nabla A|^2 + C_0 \, |A|^2 - C_1 \, t) \\ 
&\leq \max \Big \{ \sup_{x \in M_0 \cap B_4(0)} (\psi^2 \, |\nabla A|^2 + C_0 \, |A|^2),\sup_{t \in [0,T]} \sup_{x \in M_t \cap \partial B_4(0)} (C_0 \, |A|^2 - C_1 \, t) \Big \} \\ 
&\leq 1+C_0. 
\end{align*} From this, the assertion follows. \\

A similar estimate holds for the second derivatives of the second fundamental form:

\begin{lemma}
\label{shi.2}
Suppose that $M_t$, $t \in [0,T]$, is a regular mean curvature flow in $B_4(0)$ in the sense of Definition \ref{regular.flow}. Moreover, we assume that $|A| \leq 1$ for all $t \in [0,T]$ and all $x \in M_t \cap B_4(0)$. Finally, we assume that $|\nabla A|+|\nabla^2 A| \leq 1$ for all $x \in M_0 \cap B_4(0)$. Then $|\nabla A|+|\nabla^2 A| \leq C$ for all $t \in [0,1] \cap [0,T]$ and all $x \in M_t \cap B_1(0)$.
\end{lemma}

\textbf{Proof.} 
By Lemma \ref{shi.1}, we have $|\nabla A| \leq C$ for all $t \in [0,1] \cap [0,T]$ and all $x \in M_t \cap B_2(0)$. To get a bound for $|\nabla^2 A|$, we apply the maximum principle to the function $\psi^2 \, |\nabla^2 A|^2 + C_0 \, |\nabla A|^2$, where $\psi = 1-\frac{|x|^2}{4}$ and $C_0$ is a large constant. \\

Our next result will require the monotonicity formula for mean curvature flow (cf. \cite{Huisken2}). We will need a local version of this result. Specifically, we consider the modified Gaussian density 
\[\Theta(x_0,t_0;r) = \int_{M_{t_0}-r^2} \frac{1}{4\pi r^2} \, e^{-\frac{|x-x_0|^2}{4r^2}} \, (1 - |x-x_0|^2 + 4 \, r^2)_+^3.\] 
The local monotonicity formula asserts that the function $r \mapsto \Theta(x_0,t_0;r)$ is monotone increasing. A proof of this fact can be found in \cite{Ecker}, pp.~64--65 (see also \cite{Ecker-Huisken1}).

\begin{proposition} 
\label{pseudolocality.2}
There exist positive constants $\beta_0 \in (0,1)$ and $C$ such that the following holds. Suppose that $M_t$, $t \in [0,T]$, is a regular mean curvature flow in $B_4(0)$. Moreover, we assume that the initial surface $M_0$ can be expressed as the graph of a (single-valued) function $u$ over a plane. If $\|u\|_{C^4} \leq \beta_0$, then $|A(x,t)| \leq C$ for all $t \in [0,\beta_0] \cap [0,T]$ and all $x \in M_t \cap B_1(0)$.
\end{proposition}

\textbf{Proof.} 
Our argument is inspired in part by the proof of Theorem C.1 in \cite{Haslhofer-Kleiner}. Suppose that the assertion is false. Then we can find a sequence of regular mean curvature flows $\mathcal{M}_j$ in $B_4(0)$ with the following properties: 
\begin{itemize}
\item The initial surface $M_{0,j} \cap B_4(0)$ is the graph of a (single-valued) function $u_j$ over a plane, and $u_j$ satisfies $\|u_j\|_{C^4} \leq \frac{1}{j}$.
\item There exists a sequence of times $t_j \in [0,\frac{1}{j}] \cap [0,T_j]$ and a sequence of points $x_j \in M_{t_j,j} \cap B_1(0)$ such that $|A(x_j,t_j)| \geq j$.
\end{itemize}
Using a point picking argument as in Appendix C of \cite{Haslhofer-Kleiner}, we can find a pair $(\tilde{x}_j,\tilde{t}_j)$ such that $\tilde{t}_j \in [0,\frac{1}{j}] \cap [0,T_j]$, $\tilde{x}_j \in M_{\tilde{t}_j,j} \cap B_2(0)$, $Q_j := |A(\tilde{x}_j,\tilde{t}_j)| \geq j$, and 
\[\sup_{t \in [0,\tilde{t}_j]} \sup_{x \in M_{t,j} \cap B_{\frac{j}{2} \, Q_j^{-1}}(\tilde{x}_j)} |A(x,t)| \leq 2 \, Q_j.\] 
At this point, we distinguish two cases: 

\textit{Case 1:} Suppose that $\limsup_{j \to \infty} \tilde{t}_j \, Q_j^2 = 0$. By assumption, we have $|\nabla A| \leq 1$ and $|\nabla^2 A| \leq 1$ on the initial surface $M_{0,j} \cap B_4(0)$. Hence, it follows from Lemma \ref{shi.1} and Lemma \ref{shi.2} that 
\[\sup_{t \in [0,\tilde{t}_j]} \sup_{x \in M_{t,j} \cap B_{Q_j^{-1}}(\tilde{x}_j)} |\nabla A(x,t)| \leq C_0 \, Q_j^2\] 
and 
\[\sup_{t \in [0,\tilde{t}_j]} \sup_{x \in M_{t,j} \cap B_{Q_j^{-1}}(\tilde{x}_j)} |\nabla^2 A(x,t)| \leq C_0 \, Q_j^3,\] 
where $C_0$ is a uniform constant independent of $j$. In the next step, we follow the point $\tilde{x}_j$ back in time. More precisely, we consider a path $\sigma_j: [0,\tilde{t}_j] \to \mathbb{R}^3$ such that $\sigma_j(t) \in M_{t,j}$, $\sigma_j'(t)$ equals the mean curvature vector of $M_{t,j}$ at the point $\sigma_j(t)$, and $\sigma_j(\tilde{t}_j) = \tilde{x}_j$. Then $|\sigma_j'(t)| \leq 4 \, Q_j$ as long as $\sigma_j(t) \in B_{Q_j^{-1}}(\tilde{x}_j)$. Hence, if $j$ is sufficiently large, then the curve $\sigma_j(t)$ will remain in the ball $B_{Q_j^{-1}}(\tilde{x}_j)$ for all $t \in [0,\tilde{t}_j]$. In particular, if $j$ is sufficiently large, then we have $|\nabla^2 A(\sigma_j(t),t)| \leq C_0 \, Q_j^3$ for all $t \in [0,\tilde{t}_j]$. This implies 
\[\frac{d}{dt} |A(\sigma_j(t),t)| \leq C_1 \, Q_j^3\] 
for all $t \in [0,\tilde{t}_j]$, provided that $j$ is sufficiently large. Integrating this inequality from $0$ to $\tilde{t}_j$ gives 
\[Q_j = |A(\tilde{x}_j,\tilde{t}_j)| \leq |A(\sigma_j(0),0)| + C_1 \, \tilde{t}_j \, Q_j^3 \leq 1 + C_1 \, \tilde{t}_j \, Q_j^3\] 
if $j$ is sufficiently large. Since $Q_j \to \infty$ and $\tilde{t}_j \, Q_j^2 \to 0$, we arrive at a contradiction. 

\textit{Case 2:} We now assume that $\tau := \limsup_{j \to \infty} \tilde{t}_j \, Q_j^2 \in (0,\infty]$. Let us define a family of surfaces $M_{t,j}' \subset M_{t,j}$ in the following way: The surface $M_{\tilde{t}_j,j}'$ is defined as the intersection of $M_{\tilde{t}_j,j}$ with the ball $B_{\frac{j}{4} \, Q_j^{-1}}(\tilde{x}_j)$. Moreover, for each $t \in [0,\tilde{t}_j]$, the surface $M_{t,j}'$ is obtained by following each point on the surface $M_{\tilde{t}_j,j}$ back in time. It is clear that the surfaces $M_{t,j}'$, $t \in [0,\tilde{t}_j]$, form a solution of the mean curvature flow in the classical sense. Moreover, we have $\partial M_{t,j}' \cap B_{\frac{j}{8} \, Q_j^{-1}}(\tilde{x}_j) = \emptyset$ for all $t \in [0,\tilde{t}_j] \cap [\tilde{t}_j-\frac{j}{64} \, Q_j^{-2},\tilde{t}_j]$.

We next consider the rescaled surfaces $\tilde{M}_{s,j} := Q_j \, (M_{\tilde{t}_j+Q_j^{-2} \, s,j}' - \tilde{x}_j)$, $s \in [-\tilde{t}_j \, Q_j^2,0]$. These surfaces again form a solution of mean curvature flow in the classical sense. Moreover, we have $\partial \tilde{M}_{s,j} \cap B_{\frac{j}{40}}(0) = \emptyset$ for all $s \in [-\tilde{t}_j \, Q_j^2,0] \cap [-\frac{j}{64},0]$. Finally, the norm of the second fundamental form of $\tilde{M}_{s,j} \cap B_{\frac{j}{8}}(0)$ is bounded from above by $2$.

Taking the limit as $j \to \infty$, we obtain a complete, smooth, non-flat solution to the mean curvature flow which is defined on the time interval $(-\tau,0]$. The limiting solution has bounded curvature and nonnegative mean curvature. We claim that the (standard) Gaussian density of the limit flow is at most $1$. To see this, let us denote the limit flow by $\hat{M}_s$, $s \in (-\tau,0]$. 
Moreover, let us consider an arbitrary point $(y_0,s_0) \in \mathbb{R}^3 \times (-\tau,0]$ and a number $r \in (0,\sqrt{\tau+s_0})$. Clearly, $Q_j^{-1} \, r < \sqrt{\tilde{t}_j+Q_j^{-2} \, s_0}$ for $j$ large. Using Ecker's monotonicity formula for the modified Gaussian density $\Theta_{\mathcal{M}_j}$, we obtain 
\begin{align*} 
&\int_{\hat{M}_{s_0-r^2}} \frac{1}{4\pi r^2} \, e^{-\frac{|y-y_0|^2}{4r^2}} \\ 
&\leq \limsup_{j \to \infty} \Theta_{\mathcal{M}_j}(\tilde{x}_j+Q_j^{-1} \, y_0,\tilde{t}_j+Q_j^{-2} \, s_0;Q_j^{-1} \, r) \\ 
&\leq \limsup_{j \to \infty} \Theta_{\mathcal{M}_j}(\tilde{x}_j+Q_j^{-1} \, y_0,\tilde{t}_j+Q_j^{-2} \, s_0;\sqrt{\tilde{t}_j+Q_j^{-2} \, s_0}) \\ 
&= \limsup_{j \to \infty} \int_{M_{0,j}} \frac{1}{4\pi \bar{r}_j^2} \, e^{-\frac{|x-\bar{x}_j|^2}{4\bar{r}_j^2}} \, (1-|x-\bar{x}_j|^2+4\bar{r}_j^2)_+^3, 
\end{align*} 
where $\bar{x}_j := \tilde{x}_j+Q_j^{-1} \, y_0$ and $\bar{r}_j := \sqrt{\tilde{t}_j+Q_j^{-2} \, s_0}$. Using our assumption on $M_{0,j}$, we obtain 
\[\int_{\hat{M}_{s_0-r^2}} \frac{1}{4\pi r^2} \, e^{-\frac{|y-y_0|^2}{4r^2}} \leq 1\] 
for all $(y_0,s_0) \in \mathbb{R}^3 \times (-\tau,0]$ and all $r \in (0,\sqrt{\tau+s_0})$. This easily implies that the limiting solution is a flat plane of multiplicity $1$. This is a contradiction. This completes the proof of Proposition \ref{pseudolocality.2}. \\

Theorem \ref{pseudolocality} follows by combining Proposition \ref{pseudolocality.2} with Lemma \ref{shi.1} and Lemma \ref{shi.2}.

\section{The gradient estimate of Haslhofer and Kleiner (Theorem \ref{interior.derivative.estimate}
)} 

\label{haslhofer}

The curvature derivative estimate of Haslhofer and Kleiner is a consequence of the following result:

\begin{theorem}[cf. Haslhofer-Kleiner \cite{Haslhofer-Kleiner}]
\label{int.estimate}
Given any $\alpha \in (0,\frac{1}{100}]$, there exist constants $C = C(\alpha)$ and $\rho = \rho(\alpha)$ with the following property. Suppose that $M_t$, $t \in [-1,0]$, is a regular mean curvature flow in the ball $B_4(0)$. Moreover, suppose that each surface $M_t$ is outward-minimizing within the ball $B_4(0)$. We further assume that the inscribed radius and the outer radius are at least $\frac{\alpha}{H}$ at each point on $M_t$. Finally, we assume that $M_0$ passes through the origin, and $H(0,0) \leq 1$. Then $|A(x,t)| \leq C$ for all $t \in [-\rho^2,0]$ and all points $x \in M_t \cap B_\rho(0)$.
\end{theorem}

We sketch the proof of Haslhofer and Kleiner for the convenience of the reader. Suppose that there exists a sequence of regular flows $\mathcal{M}_j$ in $B_4(0)$ satisfying the assumptions of Theorem \ref{int.estimate}, and a sequence of pairs $(x_j,t_j)$ such that $t_j \in [-\frac{1}{j^2},0]$, $x_j \in M_{t_j,j} \cap B_{\frac{1}{j}}(0)$ and $|A(x_j,t_j)| \geq j^2$. By the point selection argument of Haslhofer-Kleiner, there exists a pair $(\tilde{x}_j,\tilde{t}_j)$ such that $\tilde{t}_j \in [-\frac{2}{j^2},0]$, $\tilde{x}_j \in M_{\tilde{t}_j,j} \cap B_{\frac{2}{j}}(0)$, $Q_j := |A(\tilde{x}_j,\tilde{t}_j)| \geq j^2$, and 
\[\sup_{t \in [\tilde{t}_j-\frac{j^2}{4} \, Q_j^{-2},\tilde{t}_j]} \sup_{x \in M_{t,j} \cap B_{\frac{j}{2} \, Q_j^{-1}}(\tilde{x}_j)} |A(x,t)| \leq 2 \, Q_j.\] 
We can define a family of surfaces $M_{t,j}' \subset M_{t,j}$ in the following way: The surface $M_{\tilde{t}_j,j}'$ is defined as the intersection of $M_{\tilde{t}_j,j}$ with the ball $B_{\frac{j}{4} \, Q_j^{-1}}(\tilde{x}_j)$. Moreover, for each $t \in [\tilde{t}_j-\frac{j}{64} \, Q_j^{-2},\tilde{t}_j]$, the surface $M_{t,j}'$ is obtained by following each point on the surface $M_{\tilde{t}_j,j}$ back in time. Clearly, the surfaces $M_{t,j}'$, $t \in [\tilde{t}_j-\frac{j}{64} \, Q_j^{-2},\tilde{t}_j]$, form a solution of the mean curvature flow in the classical sense. Moreover, we have $\partial M_{t,j}' \cap B_{\frac{j}{8} \, Q_j^{-1}}(\tilde{x}_j) = \emptyset$ for all $t \in [\tilde{t}_j-\frac{j}{64} \, Q_j^{-2},\tilde{t}_j]$.

Consider the rescaled surfaces $\tilde{M}_{s,j} := Q_j \, (M_{\tilde{t}_j+Q_j^{-2} \, s,j}' - \tilde{x}_j)$, $s \in [-\frac{j}{64},0]$. These surfaces again form a solution of mean curvature flow in the classical sense. Moreover, $\partial \tilde{M}_{s,j} \cap B_{\frac{j}{8}}(0) = \emptyset$ for all $s \in [-\frac{j}{64},0]$. Finally, the norm of the second fundamental form of $\tilde{M}_{s,j}$ is bounded from above by $2$. After passing to the limit as $j \to \infty$, one obtains a complete, non-flat, ancient solution to the mean curvature flow with bounded curvature. Let us denote this limit solution by $\hat{M}_s$, $s \in (-\infty,0]$.

We claim that the Gaussian density of the limit solution is at most $1$ everywhere. Let us fix a point $(y_0,s_0) \in \mathbb{R}^3 \times (-\infty,0]$ and a number $r > 0$. Using Ecker's local monotonicity formula, we obtain 
\begin{align*} 
&\int_{\hat{M}_{s_0-r^2}} \frac{1}{4\pi r^2} \, e^{-\frac{|y-y_0|^2}{4r^2}} \\ 
&\leq \limsup_{j \to \infty} \Theta_{\mathcal{M}_j}(\tilde{x}_j+Q_j^{-1} \, y_0,\tilde{t}_j+Q_j^{-2} \, s_0;Q_j^{-1} \, r) \\ 
&\leq \limsup_{j \to \infty} \Theta_{\mathcal{M}_j}(\tilde{x}_j+Q_j^{-1} \, y_0,\tilde{t}_j+Q_j^{-2} \, s_0;\sqrt{\tilde{t}_j+Q_j^{-2} \, s_0 + \frac{1}{j}}) \\ 
&= \limsup_{j \to \infty} \int_{M_{-\frac{1}{j},j}} \frac{1}{4\pi \bar{r}_j^2} \, e^{-\frac{|x-\bar{x}_j|^2}{4\bar{r}_j^2}} \, (1 - |x - \bar{x}_j|^2 + 4\bar{r}_j^2)_+^3, 
\end{align*} 
where $\bar{x}_j := \tilde{x}_j+Q_j^{-1} \, y_0$ and $\bar{r}_j := \sqrt{\tilde{t}_j+Q_j^{-2} \, s_0 + \frac{1}{j}}$.

By assumption, $\mathcal{M}_j$ satisfies $H(0,0) \leq 1$. Let $v_j$ denote the outward-pointing unit normal vector to the surface $M_{0,j}$ at the origin. Moreover, let $\Omega_{t,j}$ denote the region enclosed by $M_{t,j}$. The noncollapsing property implies that $B_\alpha(-\alpha \, v_j) \subset \Omega_{0,j} \subset \Omega_{-\frac{1}{j},j}$. On the other hand, since the surface $M_{0,j}$ passes through the origin, we must have $M_{-\frac{1}{j},j} \cap B_{\sqrt{\frac{\alpha^2}{4}+\frac{8}{j}}}(-\frac{\alpha}{2} \, v_j) \neq \emptyset$. For each $j$, we pick a point $z_j \in M_{-\frac{1}{j},j}$ which has minimal distance from $-\frac{\alpha}{2} \, v_j$ among all points on $M_{-\frac{1}{j},j}$. Then $|z_j + \frac{\alpha}{2} \, v_j| \leq \sqrt{\frac{\alpha^2}{4}+\frac{8}{j}}$ and $|z_j+\alpha \, v_j| \geq \alpha$. This implies $|z_j| \leq C \, j^{-\frac{1}{2}}$. Moreover, we have $H(z_j,-\frac{1}{j}) \leq \frac{4}{\alpha}$. In view of the noncollapsing assumption, we can find two open balls of radius $\frac{\alpha^2}{4}$ such that one of them is contained in $\Omega_{-\frac{1}{j},j}$; the other one is disjoint from $\Omega_{-\frac{1}{j},j}$; and the two balls touch each other at the point $z_j$. Consequently, the rescaled domains $j^{\frac{1}{2}} \, \Omega_{-\frac{1}{j},j}$ converge to a halfspace in the Hausdorff sense. Using the fact that $M_{-\frac{1}{j},j}$ is outward-minimizing, we conclude that the rescaled surfaces $j^{\frac{1}{2}} \, M_{-\frac{1}{j},j}$ converge, in the sense of geometric measure theory, to a plane of multiplicity at most $1$. 
%Furthermore, the outward-minimizing property implies that the area of $M_{-\frac{1}{j},j} \cap B_a(\bar{x}_j)$ is at most $4\pi a^2$ for each $a \in (0,2)$. 
Since $|\bar{x}_j| \leq O(j^{-1})$ and $\bar{r}_j = (1+o(1)) \, j^{-\frac{1}{2}}$ for $j$ large, we conclude that 
\[\limsup_{j \to \infty} \int_{M_{-\frac{1}{j},j}} \frac{1}{4\pi \bar{r}_j^2} \, e^{-\frac{|x-\bar{x}_j|^2}{4\bar{r}_j^2}} \, (1 - |x - \bar{x}_j|^2 + 4\bar{r}_j^2)_+^3 \leq 1.\] 
Therefore, the limiting flow $\hat{M}_s$ has Gaussian density at most $1$. This contradicts the fact that the limit flow is non-flat.

\section{Proof of Theorem \ref{properties.of.surgery}}

\label{construction.of.cap}

In this section, we explain our procedure for capping off a neck. We begin by constructing an axially symmetric model surface.

\begin{lemma}
\label{model.surface}
The surface 
\[\Sigma = \Big \{ \Big ( \sqrt{\frac{s}{1+s}} \, \cos(2\pi t),\sqrt{\frac{s}{1+s}} \, \sin(2\pi t),s \Big ) : s \in [0,\infty), \, t \in [0,1] \Big \}\] 
closes up smoothly at $s=0$. Moreover, we have $0 < \lambda_1 < \lambda_2 = \mu$ whenever $s>0$. Here, $\mu$ denotes the reciprocal of the inscribed radius of $\Sigma$.
\end{lemma}

\textbf{Proof.} 
The smoothness of $\Sigma$ is obvious. A straightforward calculation shows that the principal curvatures of $\Sigma$ are given by 
\[\lambda_1 = 2 \, (1+4s \, (1+s)^3)^{-\frac{3}{2}} \, (1+s)^2 \, (1+4s)\] 
and 
\[\lambda_2 = 2 \, (1+4s \, (1+s)^3)^{-\frac{1}{2}} \, (1+s)^2.\] 
Clearly, $0 < \lambda_1 < \lambda_2$ for $s > 0$. Hence, it remains to estimate the inscribed radius of $\Sigma$. To that end, let $U = \{x \in \mathbb{R}^3: x_3 > 0, \, x_1^2+x_2^2 < \frac{x_3}{1+x_3}\}$ be the region enclosed by $\Sigma$. For each $s>0$, we denote by $W_s$ the open ball of radius $\frac{1}{2} \, (1+4s \, (s+1)^3)^{1/2} \, (1+s)^{-2}$ centered at the point $(0,0,s+\frac{1}{2} \, (s+1)^{-2})$. For each $s>0$, the circle 
\[C_s := \Big \{ \Big ( \sqrt{\frac{s}{1+s}} \, \cos(2\pi t),\sqrt{\frac{s}{1+s}} \, \sin(2\pi t),s \Big ) : t \in [0,1] \Big \}\] 
is contained in $\Sigma \cap \partial W_s$. Moreover, the surfaces $\Sigma$ and $\partial W_s$ have the same tangent plane at each point on the circle $C_s$.

It is easy to see that $W_s \subset U$ if $s$ is sufficiently large. We claim that $W_s \subset U$ for all $s > 0$. Suppose this is false. Let $\bar{s} = \sup \{s > 0: W_s \not\subset U\}$. Then $W_{\bar{s}} \subset U$. Moreover, we can find a sequence of numbers $s_j \nearrow \bar{s}$ and a sequence of points $p_j \in W_{s_j} \setminus U$. After passing to a subsequence if necessary, the points $p_j$ converge to some point $p \in \bar{W}_{\bar{s}} \setminus U$. Since $W_{\bar{s}} \subset U$, we conclude that $p \in \Sigma \cap \partial W_{\bar{s}}$, and the surfaces $\Sigma$ and $\partial W_{\bar{s}}$ have the same tangent plane at the point $p$. On the other hand, since $\lambda_1 < \lambda_2$, we must have $\liminf_{j \to \infty} \text{\rm dist}(p_j,C_{s_j}) > 0$. Consequently, we have $p \in C_{\tilde{s}}$ for some $\tilde{s} \neq \bar{s}$. This implies that $p \in \Sigma \cap \partial W_{\tilde{s}}$, and the surfaces $\Sigma$ and $\partial W_{\tilde{s}}$ have the same tangent plane at the point $p$. Thus, the spheres $\partial W_{\bar{s}}$ and $\partial W_{\tilde{s}}$ touch each other at the point $p$, but this is impossible if $\bar{s} \neq \tilde{s}$. This shows that $W_s \subset U$ for all $s > 0$. Consequently, the inscribed radius is given by $\frac{1}{2} \, (1+4s \, (s+1)^3)^{1/2} \, (1+s)^{-2}$. This completes the proof of Lemma \ref{model.surface}. \\

In the remainder of this section, we consider an $(\hat{\alpha},\hat{\delta},\varepsilon,L)$-neck $N$ of size $1$, which is contained in a closed, embedded, mean convex surface $M \subset \mathbb{R}^3$. It is understood that $\varepsilon$ is much smaller than $\hat{\delta}$. By definition, we can find a simple closed, convex curve $\Gamma$ with the property that $\text{\rm dist}_{C^{20}}(N,\Gamma \times [-L,L]) \leq \varepsilon$. Moreover, the curve $\Gamma$ is $\frac{1}{1+\hat{\delta}}$-noncollapsed, and the derivatives of the geodesic curvature of $\Gamma$ satisfy $\sum_{l=1}^{18} |\nabla^l \kappa| \leq \frac{1}{100}$ at each point on $\Gamma$. Furthermore, there exists a point on $\Gamma$ where the geodesic curvature $\kappa$ is equal to $1$.

Since $\text{\rm dist}_{C^{20}}(N,\Gamma \times [-L,L]) \leq \varepsilon$, we can find a collection of curves $\Gamma_s$ such that 
\[\{(\gamma_s(t),s): s \in [-(L-1),L-1], \, t \in [0,1]\} \subset N\] 
and 
\[\sum_{k+l \leq 20} \Big | \frac{\partial^k}{\partial s^k} \, \frac{\partial^l}{\partial t^l} (\gamma_s(t) - \gamma(t)) \Big | \leq O(\varepsilon).\]
Here, we have used the notation $\Gamma = \{\gamma(t): t \in [0,1]\}$ and $\Gamma_s = \{\gamma_s(t): t \in [0,1]\}$.

The following lemma is analogous to Proposition 3.17 in \cite{Huisken-Sinestrari3}: 

\begin{lemma}[cf. Huisken-Sinestrari \cite{Huisken-Sinestrari3}, Proposition 3.17]
\label{aux}
Consider a bended surface of the form 
\[\tilde{N} = \{((1-u(s)) \, \gamma_s(t),s): s \in (0,\Lambda^{\frac{1}{4}}], \, t \in [0,1]\},\] 
where $|u|+|u'|+|u''| \leq \frac{1}{10}$ everywhere.
Then we have the pointwise estimates 
\[\tilde{\lambda}_1(s,t) \geq \lambda_1(s,t) + c_0 \, u''(s) - \frac{1}{c_0} \, (|u(s)|+|u'(s)|)\] 
and 
\[\tilde{H}(s,t) \geq H(s,t) + c_0 \, u''(s) - \frac{1}{c_0} \, (|u(s)|+|u'(s)|),\] 
where $c_0>0$ is a universal constant.
\end{lemma}

It will be convenient to translate the neck $N$ in space so that the center of mass of $\Gamma$ is at the origin. Using the curve shortening flow, we can construct a homotopy $\tilde{\gamma}_r(t)$, $(r,t) \in [0,1] \times [0,1]$, with the following properties: 
\begin{itemize} 
\item $\tilde{\gamma}_r(t)=\gamma(t)$ for $r \in [0,\frac{1}{4}]$.
\item $\tilde{\gamma}_r(t)=(\cos (2\pi t),\sin(2\pi t))$ for $r \in [\frac{1}{2},1]$.
\item For each $r \in [0,1]$, the curve $\tilde{\Gamma}_r$ is $\frac{1}{1+\hat{\delta}}$-noncollapsed.
\item We have $\sup_{(r,t) \in [0,1] \times [0,1]} |\frac{\partial}{\partial r} \tilde{\gamma}_r(t)| + |\frac{\partial^2}{\partial r \, \partial t} \tilde{\gamma}_r(t)| + |\frac{\partial^2}{\partial r^2} \tilde{\gamma}_r(t)| \leq \omega(\hat{\delta})$, where $\omega(\hat{\delta}) \to 0$ as $\hat{\delta} \to 0$.
\end{itemize} 
Moreover, let $\chi: \mathbb{R} \to \mathbb{R}$ be a smooth cutoff function such that $\chi = 1$ on $(-\infty,1]$ and $\chi=0$ on $[2,\infty)$. We next define a surface $\tilde{F}_\Lambda: [-L,\Lambda] \times [0,1] \to \mathbb{R}^3$ by 
\[\tilde{F}_\Lambda(s,t) = \begin{cases} (\gamma_s(t),s) & \text{\rm for $s \in [-(L-1),0]$} \\ ((1-e^{-\frac{4\Lambda}{s}}) \, \gamma_s(t),s) & \text{\rm for $s \in (0,\Lambda^{\frac{1}{4}}]$} \\ ((1-e^{-\frac{4\Lambda}{s}}) \, (\chi(s/\Lambda^{\frac{1}{4}}) \, \gamma_s(t)+(1-\chi(s/\Lambda^{\frac{1}{4}})) \, \gamma(t)),s) & \text{\rm for $s \in (\Lambda^{\frac{1}{4}},2 \, \Lambda^{\frac{1}{4}}]$} \\ ((1-e^{-\frac{4\Lambda}{s}}) \, \tilde{\gamma}_{s/\Lambda}(t),s) & \text{\rm for $s \in (2 \, \Lambda^{\frac{1}{4}},\Lambda]$.} \end{cases}\] 
It is clear that $\tilde{F}_\Lambda$ is smooth. Moreover, $\tilde{F}_\Lambda$ is axially symmetric for $s \geq \frac{\Lambda}{2}$. \\

\begin{lemma}
\label{convexity.of.cap}
We can find real numbers numbers $\delta_1$, $\Lambda_1$, and a function $E(\Lambda)$ such that the following statements hold: 
\begin{itemize} 
\item Suppose that $\hat{\delta} < \delta_1$, $\frac{L}{1000} \geq \Lambda \geq \Lambda_1$, and $\varepsilon \leq E(\Lambda)$. Then, for each point $(s,t) \in (0,\Lambda^{\frac{1}{4}}] \times [0,1]$, the mean curvature of $\tilde{F}_\Lambda$ at $(s,t)$ is greater than the mean curvature of the original neck at $(s,t)$, and the smallest curvature eigenvalue of $\tilde{F}_\Lambda$ is greater than the smallest curvature eigenvalue of the original neck at $(s,t)$. 
\item Suppose that $\hat{\delta} < \delta_1$, $\frac{L}{1000} \geq \Lambda \geq \Lambda_1$, and $\varepsilon \leq E(\Lambda)$. Then, for each point $(s,t) \in (\Lambda^{\frac{1}{4}},2 \, \Lambda^{\frac{1}{4}}] \times [0,1]$, the surface $\tilde{F}_\Lambda$ is strictly convex at $(s,t)$.
\item Suppose that $\hat{\delta} < \delta_1$, $\frac{L}{1000} \geq \Lambda \geq \Lambda_1$, and $\varepsilon \leq E(\Lambda)$. Then, for each point $(s,t) \in (2 \, \Lambda^{\frac{1}{4}},\Lambda] \times [0,1]$, the surface $\tilde{F}_\Lambda$ is strictly convex at $(s,t)$.
\end{itemize}
\end{lemma}

\textbf{Proof.} 
We begin with the first statement. Let $\tilde{\lambda}_1$ denote the smallest curvature eigenvalue of the bended surface $\tilde{F}_\Lambda$ and let $\lambda_1$ be the smallest curvature eigenvalue of the original neck. Similarly, we denote by $\tilde{H}$ the mean curvature of the bended surface and by $H$ the mean curvature of the original neck. By choosing $\Lambda$ sufficiently large, we can arrange that the function $u(s) = e^{-\frac{4\Lambda}{s}}$ satisfies 
\[u''(s) \geq \frac{1}{c_0^2} \, (|u(s)|+|u'(s)|)\] 
for all $s \in (0,\Lambda^{\frac{1}{4}}]$, where $c_0$ is the constant from Lemma \ref{aux}. Using Lemma \ref{aux}, we conclude that $\tilde{\lambda}_1 \geq \lambda_1$ and $\tilde{H} \geq H$ for all points $(s,t) \in (0,\Lambda^{\frac{1}{4}}] \times [0,1]$. This proves the first statement.

To verify the second statement, we consider a point $(s,t) \in (\Lambda^{\frac{1}{4}},2 \, \Lambda^{\frac{1}{4}}] \times [0,1]$. It is easy to see that $\tilde{h}_{tt} \geq \frac{1}{2}$ and $\langle (\gamma(t),0),\tilde{\nu}(s,t) \rangle \geq \frac{1}{2}$. We next compute 
\[\frac{\partial^2}{\partial s^2} \tilde{F}_\Lambda(s,t) = -\Big ( \frac{16\Lambda^2}{s^4} - \frac{8\Lambda}{s^3} \Big ) \, e^{-\frac{4\Lambda}{s}} \, (\gamma(t),0) + O(\varepsilon)\] 
and 
\begin{align*} 
\frac{\partial^2}{\partial s \, \partial t} \tilde{F}_\Lambda(s,t) 
&= -\frac{4\Lambda}{s^2} \, e^{-\frac{4\Lambda}{s}} \, \Big ( \chi(s/\Lambda^{\frac{1}{4}}) \, \frac{\partial}{\partial t} \gamma_s(t) + (1-\chi(s/\Lambda^{\frac{1}{4}})) \, \frac{\partial}{\partial t} \gamma(t),0 \Big ) + O(\varepsilon) \\ 
&= -\frac{4\Lambda}{s^2} \, \frac{e^{-\frac{4\Lambda}{s}}}{1-e^{-\frac{4\Lambda}{s}}} \, \frac{\partial \tilde{F}_\Lambda}{\partial t}(s,t) + O(\varepsilon) 
\end{align*} 
for $(s,t) \in (\Lambda^{\frac{1}{4}},2 \, \Lambda^{\frac{1}{4}}] \times [0,1]$. From this, we deduce that 
\begin{align*} 
\tilde{h}_{ss} 
&= -\Big \langle \frac{\partial^2}{\partial s^2} \tilde{F}_\Lambda(s,t),\tilde{\nu}(s,t) \Big \rangle \\
&= \Big ( \frac{16\Lambda^2}{s^4} - \frac{8\Lambda}{s^3} \Big ) \, e^{-\frac{4\Lambda}{s}} \, \langle (\gamma(t),0),\tilde{\nu}(s,t) \rangle + O(\varepsilon) \\ 
&\geq \Big ( \frac{8\Lambda^2}{s^4} - \frac{4\Lambda}{s^3} \Big ) \, e^{-\frac{4\Lambda}{s}} + O(\varepsilon) 
\end{align*}
and 
\[\tilde{h}_{st} = -\Big \langle \frac{\partial^2}{\partial s \, \partial t} \tilde{F}_\Lambda(s,t),\tilde{\nu}(s,t) \Big \rangle = O(\varepsilon)\] 
for $(s,t) \in (\Lambda^{\frac{1}{4}},2 \, \Lambda^{\frac{1}{4}}] \times [0,1]$. Hence, if $\varepsilon$ is small enough (depending on $\Lambda$), then $\tilde{h}_{ss} \, \tilde{h}_{tt} - \tilde{h}_{st}^2 > 0$, and the surface $\tilde{F}_\Lambda$ is strictly convex at $(s,t)$. This completes the proof of the second statement.

To prove the third statement, we consider a point $(s,t) \in (2 \, \Lambda^{\frac{1}{4}},\Lambda] \times [0,1]$. We clearly have $\tilde{h}_{tt} \geq \frac{1}{2}$ and $\langle (\tilde{\gamma}_{s/\Lambda}(t),0),\tilde{\nu}(s,t) \rangle \geq \frac{1}{2}$. We next compute 
\[\frac{\partial^2}{\partial s^2} \tilde{F}_\Lambda(s,t) = -\Big ( \frac{16\Lambda^2}{s^4} - \frac{8\Lambda}{s^3} \Big ) \, e^{-\frac{4\Lambda}{s}} \, (\tilde{\gamma}_{s/\Lambda}(t),0) + O(\Lambda^{-2} \, \omega(\hat{\delta}) \, 1_{\{\frac{\Lambda}{4} \leq s \leq \frac{\Lambda}{2}\}})\] 
and 
\begin{align*} 
\frac{\partial^2}{\partial s \, \partial t} \tilde{F}_\Lambda(s,t) 
&= -\frac{4\Lambda}{s^2} \, e^{-\frac{4\Lambda}{s}} \, \Big ( \frac{\partial}{\partial t} \tilde{\gamma}_{s/\Lambda}(t),0 \Big ) + O(\Lambda^{-1} \, \omega(\hat{\delta}) \, 1_{\{\frac{\Lambda}{4} \leq s \leq \frac{\Lambda}{2}\}}) \\ 
&= -\frac{4\Lambda}{s^2} \, \frac{e^{-\frac{4\Lambda}{s}}}{1-e^{-\frac{4\Lambda}{s}}} \, \frac{\partial \tilde{F}_\Lambda}{\partial t}(s,t) + O(\Lambda^{-1} \, \omega(\hat{\delta}) \, 1_{\{\frac{\Lambda}{4} \leq s \leq \frac{\Lambda}{2}\}}) 
\end{align*} 
for $(s,t) \in (2 \, \Lambda^{\frac{1}{4}},\Lambda] \times [0,1]$. From this, we deduce that 
\begin{align*} 
\tilde{h}_{ss} 
&= -\Big \langle \frac{\partial^2}{\partial s^2} \tilde{F}_\Lambda(s,t),\tilde{\nu}(s,t) \Big \rangle \\ 
&= \Big ( \frac{16\Lambda^2}{s^4} - \frac{8\Lambda}{s^3} \Big ) \, e^{-\frac{4\Lambda}{s}} \, \langle (\tilde{\gamma}_{s/\Lambda}(t),0),\tilde{\nu}(s,t) \rangle + O(\Lambda^{-2} \, \omega(\hat{\delta}) \, 1_{\{\frac{\Lambda}{4} \leq s \leq \frac{\Lambda}{2}\}}) \\ 
&\geq \Big ( \frac{8\Lambda^2}{s^4} - \frac{4\Lambda}{s^3} \Big ) \, e^{-\frac{4\Lambda}{s}} + O(\Lambda^{-2} \, \omega(\hat{\delta}) \, 1_{\{\frac{\Lambda}{4} \leq s \leq \frac{\Lambda}{2}\}}) 
\end{align*}
and 
\[\tilde{h}_{st} = -\Big \langle \frac{\partial^2}{\partial s \, \partial t} \tilde{F}_\Lambda(s,t),\tilde{\nu}(s,t) \Big \rangle = O(\Lambda^{-1} \, \omega(\hat{\delta}) \, 1_{\{\frac{\Lambda}{4} \leq s \leq \frac{\Lambda}{2}\}})\] 
for $(s,t) \in (2 \, \Lambda^{\frac{1}{4}},\Lambda] \times [0,1]$. Hence, if $\hat{\delta}$ is sufficiently small, then $\tilde{h}_{ss} \, \tilde{h}_{tt} - \tilde{h}_{st}^2 > 0$, and the surface $\tilde{F}_\Lambda$ is strictly convex at $(s,t)$. This completes the proof of Lemma \ref{convexity.of.cap}. \\

Since the surface $\tilde{F}_\Lambda$ is axially symmetric in the region $\{\frac{\Lambda}{2} \leq s \leq \Lambda\}$, we may glue $\tilde{F}_\Lambda$ to a scaled copy of the axially symmetric cap constructed in Lemma \ref{model.surface}. We briefly sketch how this can be done. Let us fix a smooth, convex, even function $\Phi: \mathbb{R} \to \mathbb{R}$ such that $\Phi(z) = |z|$ for $|z| \geq \frac{1}{100}$. For $\Lambda$ very large, we define 
\[a = 1-e^{-4}+\frac{1}{3} \, (1-e^{-4})^2 \, \Lambda^{-\frac{1}{4}}\]
and 
\begin{align*} 
v_\Lambda(s) 
&= 1-e^{-\frac{4\Lambda}{s}} + a \, \sqrt{\frac{\Lambda+2\,\Lambda^{\frac{1}{4}}-s}{a+\Lambda+2\,\Lambda^{\frac{1}{4}}-s}} \\ 
&- \Lambda^{-\frac{1}{4}} \, \Phi \bigg ( \Lambda^{\frac{1}{4}} \, \Big ( 1-e^{-\frac{4\Lambda}{s}} - a \, \sqrt{\frac{\Lambda+2\,\Lambda^{\frac{1}{4}}-s}{a+\Lambda+2\,\Lambda^{\frac{1}{4}}-s}} \Big ) \bigg ) 
\end{align*} 
for $s \in [\Lambda,\Lambda+\Lambda^{\frac{1}{4}}]$. Since the functions $s \mapsto 1-e^{-\frac{4\Lambda}{s}}$ and $s \mapsto a \, \sqrt{\frac{\Lambda+2 \, \Lambda^{\frac{1}{4}}-s}{a+\Lambda+2\,\Lambda^{\frac{1}{4}}-s}}$ are concave, the function $v_\Lambda$ is concave as well. Moreover, if $\Lambda$ is sufficiently large, then we have $v_\Lambda(s)=2 \, (1-e^{-\frac{4\Lambda}{s}})$ in a neighborhood of the point $s=\Lambda$, and $v_\Lambda(s)=2a \, \sqrt{\frac{\Lambda+2\,\Lambda^{\frac{1}{4}}-s}{a+\Lambda+2\,\Lambda^{\frac{1}{4}}-s}}$ in a neighborhood of the point $s=\Lambda+\Lambda^{\frac{1}{4}}$. We now extend $\tilde{F}_\Lambda$ to the region $[-(L-1),\Lambda+2\Lambda^{\frac{1}{4}}] \times [0,1]$ by putting 
\[\tilde{F}_\Lambda(s,t) = \Big ( \frac{1}{2} \, v_\Lambda(s) \, \cos(2\pi t),\frac{1}{2} \, v_\Lambda(s) \, \sin(2\pi t),s \Big )\] 
for $s \in (\Lambda,\Lambda+\Lambda^{\frac{1}{4}}]$ and 
\[\tilde{F}_\Lambda(s,t) = \Big ( a \, \sqrt{\frac{\Lambda+2 \, \Lambda^{\frac{1}{4}} - s}{a+\Lambda+2\,\Lambda^{\frac{1}{4}}-s}} \, \cos(2\pi t),a \, \sqrt{\frac{\Lambda+2 \, \Lambda^{\frac{1}{4}} - s}{a+\Lambda+2\,\Lambda^{\frac{1}{4}}-s}} \, \sin(2\pi t),s \Big )\] 
for $s \in (\Lambda+\Lambda^{\frac{1}{4}},\Lambda+2\Lambda^{\frac{1}{4}}]$. It is straightforward to verify that the resulting surface is smooth satisfies the curvature bounds $\frac{1}{2} \leq H \leq 10$. 

\begin{lemma}
The surface $\tilde{F}_\Lambda$ is convex in the region $(\Lambda,\Lambda+2\Lambda^{\frac{1}{4}}] \times [0,1]$.
\end{lemma}

\textbf{Proof.} 
Since the function $v_\Lambda$ is concave on the interval $[\Lambda,\Lambda+\Lambda^{\frac{1}{4}}]$, we conclude that $\tilde{F}_\Lambda$ is a convex surface for $(s,t) \in (\Lambda,\Lambda+\Lambda^{\frac{1}{4}}] \times [0,1]$. Moreover, it follows from Lemma \ref{model.surface} that the surface $\tilde{F}_\Lambda$ is convex for $(s,t) \in (\Lambda+\Lambda^{\frac{1}{4}},\Lambda+2\Lambda^{\frac{1}{4}}]$. \\

In the next step, we show that in the surgically modified region the inscribed radius is at least $\frac{1}{(1+\delta) \, H}$ and the outer radius is at least $\frac{\alpha}{H}$.

\begin{proposition}
\label{noncollapsing.preserved.under.surgery}
Given any number $\hat{\alpha} > \alpha$, we can find a number $\delta_2$ with the following property. Suppose that we are given a pair of real numbers $\delta$ and $\hat{\delta}$ such that $\hat{\delta} < \delta < \delta_2$. Then there exist real numbers $\bar{\varepsilon}$ and $\Lambda_2$ such that the surface $\tilde{F}_\Lambda$ is $\frac{1}{1+\delta}$-noncollapsed whenever $\varepsilon \leq \bar{\varepsilon}$ and $\frac{L}{1000} \geq \Lambda \geq \Lambda_2$. Furthermore, if $\varepsilon \leq \bar{\varepsilon}$ and $\frac{L}{1000} \geq \Lambda \geq \Lambda_2$, then the outer radius is at least $\frac{\alpha}{H}$ at each point on $\tilde{F}_\Lambda$.
\end{proposition}

\textbf{Proof.} 
We first establish the bound for the inscribed radius. It follows from Lemma \ref{model.surface} that the inscribed radius of $\tilde{F}_\Lambda$ is at least $\frac{1}{(1+\delta) \, H}$ at each point in the region $(\Lambda+\Lambda^{\frac{1}{4}},\Lambda+2\Lambda^{\frac{1}{4}}]$. Consider now a point $(s_0,t_0) \in [-(L-1),\Lambda+\Lambda^{\frac{1}{4}}] \times [0,1]$. If $\Lambda$ is large, then we can approximate the map $\tilde{F}_\Lambda$ near the point $(s_0,t_0)$ by a cylinder whose cross-section is a simple closed, convex curve. Furthermore, the noncollapsing constant of the cross section is $\frac{1}{1+\hat{\delta}}-O(\varepsilon)$ or better. Since $\hat{\delta} < \delta$, we conclude that the surface $\tilde{F}_\Lambda$ is $\frac{1}{1+\delta}$-noncollapsed if $\Lambda$ is sufficiently large and $\varepsilon$ is sufficiently small. 

It remains to prove the bound for the outer radius. Let $\tilde{N}$ denote the image of the map $\tilde{F}_\Lambda: (-(,\Lambda+2\Lambda^{\frac{1}{4}}] \times [0,1] \to \mathbb{R}^3$. By assumption, the region $\{x + a \, \nu(x): x \in N, \, a \in (0,2\hat{\alpha})\}$ is disjoint from $M \setminus N$. Since the surface $\tilde{N} \setminus N$ lies inside the original neck $N$, it follows that the region $\{x + a \, \tilde{\nu}(x): x \in \tilde{N} \setminus N, \, a \in (0,2\hat{\alpha})\}$ is disjoint from $M \setminus N$. Consequently, for each point $x \in \tilde{N} \setminus N$, we can find a ball of radius $\hat{\alpha}$ which touches $\tilde{N}$ at the point $x$ from the outside, and which is disjoint from $M \setminus N$. On the other hand, if $\hat{\delta}$ and $\varepsilon$ are sufficiently small and $\Lambda$ is sufficiently large, then the mean curvature of the surface $\tilde{N} \setminus N$ is greater than $\frac{\alpha}{\hat{\alpha}}$ everywhere. Putting these facts together, we conclude that the outer radius is at least $\hat{\alpha} > \frac{\alpha}{H}$ at each point in $\tilde{N} \setminus N$. This completes the proof of Proposition \ref{noncollapsing.preserved.under.surgery}. \\

We note that our surgery procedure always produces an embedded surface. Finally, it is clear from the construction that the resulting cap is at least of class $C^5$ with uniform bounds independent of the surgery parameters $\hat{\alpha}$, $\hat{\delta}$, $\varepsilon$, $L$, and $\Lambda$.

\section{Proof of Proposition \ref{separation.of.surgery.regions}}

By assumption, the point $x_0$ lies in the surgically modified region of $M_{t_0+}$. Hence, the surface $M_{t_0-}$ must have contained an $(\hat{\alpha},\hat{\delta},\varepsilon,L)$-neck of size $r \in [\frac{1}{2H_1},\frac{2}{H_1}]$. Let us denote this neck by $N$. At time $t_0$ the neck $N$ is  replaced by a capped-off neck $\tilde{N}$. More precisely, suppose that the original neck $N$ satisfies 
\[\{(\gamma_s(t),s): s \in [-(L-1),L-1], \, t \in [0,1]\} \subset r^{-1} \, N.\] 
Then the surface $\tilde{N}$ satisfies 
\[N \cap \{x \in \mathbb{R}^3: \langle x,e_3 \rangle \leq 0\} \subset \tilde{N}\] 
and 
\[\tilde{N} \subset \{x \in \mathbb{R}^3: \langle x,e_3 \rangle \leq 4\Lambda \, H_1^{-1}\}.\] 
Since the point $x_0$ lies in the surgically modified part of $\tilde{N}$, we have $\langle x_0,e_3 \rangle \geq 0$.

By assumption, the outer radius of $M_{t_0+}$ is at least $\frac{\alpha}{H}$ everywhere. Moreover, we have $H \leq 100 \, H_1$ at each point on $\tilde{N}$. Therefore, for each point $x \in \tilde{N}$, the outer radius is at least $\frac{\alpha}{100} \, H_1^{-1}$. Hence, if we denote by $\nu$ the outward-pointing unit normal vector field to $\tilde{N}$, then the set  
\[E = \{x + a \, \nu(x): x \in \tilde{N}, \, a \in (0,\frac{\alpha}{100} \, H_1^{-1})\}\] 
is disjoint from the region $\Omega_{t_0+}$. 

By assumption, the point $x_1 \in M_{t_1+}$ lies in the surgically modified region at time $t_1$. This region was created by performing surgery on an $(\hat{\alpha},\hat{\delta},\varepsilon,L)$-neck in $M_{t_1-}$. This neck has length at least $\frac{L}{2} \, H_1^{-1} \geq 50 \, \Lambda \, H_1^{-1}$. Hence, we can find two points $y,z \in \mathbb{R}^3$ such that $|y-z| = 40 \, \Lambda \, H_1^{-1}$, $|\frac{y+z}{2}-x_1| \leq \frac{\alpha}{1000} \, H_1^{-1}$, and the line segment joining $y$ and $z$ is contained in the region $\Omega_{t_1-}$. Since $\Omega_{t_1-} \subset \Omega_{t_0+}$ is disjoint from $E$, the line segment joining $y$ and $z$ cannot intersect the set $E$. 

Now, if $|\frac{y+z}{2}-x_0| \leq \frac{\alpha}{500} \, H_1^{-1}$, then it is easy to see that the line segment joining $y$ and $z$ must intersect the set $E$. Therefore, we have $|\frac{y+z}{2}-x_0| > \frac{\alpha}{500} \, H_1^{-1}$, hence $|x_1-x_0| > \frac{\alpha}{1000} \, H_1^{-1}$. This completes the proof of Proposition \ref{separation.of.surgery.regions}.

\section{Proof of Proposition \ref{consequence.of.pseudolocality}}

Let $\beta_0$ be chosen as in Theorem \ref{pseudolocality}. Let us fix a positive number $\beta_1 \in (0,\frac{\alpha}{8000})$ such that the following holds: Suppose that $N$ is an $(\hat{\alpha},\hat{\delta},\varepsilon,L)$-neck of size $r \in [\frac{1}{2H_1},\frac{2}{H_1}]$. Moreover, let $\tilde{N}$ denote the capped-off neck obtained by performing a $\Lambda$-surgery on $N$. Then, for each point $x_0$ in the surgically modified region, the dilated surface $(\beta_1^{-1} \, H_1 \, (\tilde{N} - x_0)) \cap B_4(0)$ can be expressed as the graph of function which has $C^4$-norm less than $\beta_0$. In view of the construction of the cap in Section \ref{construction.of.cap}, we can choose the constant $\beta_1$ in such a way that $\beta_1$ depends only on the noncollapsing constant $\alpha$, but not on the exact choice of the surgery parameters $\hat{\alpha}$, $\hat{\delta}$, $\varepsilon$, $L$, and $H_1$.

After these preparations, we now complete the proof of Proposition \ref{consequence.of.pseudolocality}. Suppose that $t_0$ is a surgery time and $x_0$ lies in the surgically modified region. By Proposition \ref{separation.of.surgery.regions}, the flow $M_t \cap B_{4\beta_1 \, H_1^{-1}}(x_0)$ is smooth for all times $t > t_0$. Moreover, the surface $(\beta_1^{-1} \, H_1 \, (M_{t_0+} - x_0)) \cap B_4(0)$ is a graph of a function with $C^4$-norm less than $\beta_0$. Hence, Theorem \ref{pseudolocality} implies that 
\[\beta_1 \, H_1^{-1} \, |A| + \beta_1^2 \, H_1^{-2} \, |\nabla A| + \beta_1^3 \, H_1^{-3} \, |\nabla^2 A| \leq C\] 
for all $t \in (t_0,t_0+\beta_0 \, \beta_1^2 \, H_1^{-2}]$ and all $x \in M_t \cap B_{\beta_1 \, H_1^{-1}}(x_0)$. From this, the assertion follows.

\section{Proof of Proposition \ref{gradient.estimate}}

Let us consider an arbitrary time $t_1 \geq (1000 \, \sup_{M_0} |A|)^{-2}$ and an arbitrary point $x_1 \in M_{t_1}$ for which we want to verify the estimate. To avoid confusion (and without any loss of generality), we will assume that $t_1$ is not itself a surgery time. There are two cases: 

\textit{Case 1:} There exists a surgery time $t_0$ and a point $x_0$ such that $|x_1-x_0| \leq \beta_* \, H_1^{-1}$, $0 < t_1-t_0 \leq \beta_* \, H_1^{-2}$, and $x_0$ lies in the surgically modified region at time $t_0+$. Applying Proposition \ref{consequence.of.pseudolocality}, we conclude that 
\[H_1^{-1} \, |A| + H_1^{-2} \, |\nabla A| + H_1^{-3} \, |\nabla^2 A| \leq C_*\] 
at the point $(x_1,t_1)$. Hence, $|\nabla A| \leq C_* \, (H+H_1)^2$ and $|\nabla^2 A| \leq C_* \, (H+H_1)^3$ at the point $(x_1,t_1)$. 

\textit{Case 2:} There does not exist a surgery time $t_0$ and a point $x_0$ such that $|x_1-x_0| \leq \beta_* \, H_1^{-1}$, $0 < t_1-t_0 \leq \beta_* \, H_1^{-2}$, and $x_0$ lies in the surgically modified region at time $t_0+$. In this case, the surfaces $M_t \cap B_{\beta_* \, H_1^{-1}}(x_1)$, $t \in (t_1-\beta_* \, H_1^{-2},t_1]$, form a regular mean curvature flow in the sense of Definition \ref{regular.flow}. Note that, since $t_1 \geq (1000 \, \sup_{M_0} |A|)^{-2}$, we have $t_1-\beta_* \, H_1^{-2} > 0$. Moreover, the ball $B_{\beta_* \, H_1^{-1}}(x_1)$ is contained in the region $\Omega_0$, so the surfaces $M_t$ are outward-minimizing within the ball $B_{\beta_* \, H_1^{-1}}(x_1)$. Hence, Theorem \ref{interior.derivative.estimate} implies that $|\nabla A| \leq B \, (H+H_1)^2$ and $|\nabla^2 A| \leq B \, (H+H_1)^3$ at the point $(x_1,t_1)$. Here, $B$ is a positive constant that depends only on $\beta_*$ and the noncollapsing constant $\alpha$. This completes the proof.

\section{Proof of Proposition \ref{choice.of.delta}}

We next prove some auxiliary results about curves. In the following, we assume that $C_\#$ is the constant in Proposition \ref{gradient.estimate}.

\begin{lemma}
\label{const.curvature}
Let $\Gamma$ be a (possibly non-closed) embedded curve in the plane of class $C^3$ with geodesic curvature $\kappa>0$. Moreover, suppose that the inscribed radius is at least $\frac{1}{\kappa}$ at each point on $\Gamma$. Then $\kappa$ is constant.
\end{lemma}

\textbf{Proof.} 
The assumption implies that the function 
\[Z(s,t) := \frac{1}{2} \, \kappa(s) \, |\gamma(s)-\gamma(t)|^2 - \langle \gamma(s)-\gamma(t),\nu(s) \rangle\] 
is nonnegative for all $s,t$. A straightforward calculation gives 
\[\frac{\partial Z}{\partial t}(s,t) \Big |_{s=t} = 0, \quad \frac{\partial^2 Z}{\partial t^2}(s,t) \Big |_{s=t} = 0, \quad \frac{\partial^3 Z}{\partial t^3}(s,t) \Big |_{s=t} = -\frac{d\kappa}{ds}(s).\] 
Since $Z$ is nonnegative everywhere, we conclude that $\frac{d\kappa}{ds}(s) = 0$ at each point on $\Gamma$.

\begin{lemma} 
\label{sup}
Let $\Gamma_j$ be a sequence of (possibly non-closed) embedded curves in the plane with the property that $\kappa > 0$, $|\frac{d\kappa}{ds}| \leq C_\# \, (\kappa+2\Theta)^2$, and $|\frac{d^2\kappa}{ds^2}| \leq C_\# \, (\kappa+2\Theta)^3$. Moreover, suppose that the inscribed radius is at least $\frac{1}{(1+\frac{1}{j}) \, \kappa}$ at each point on $\Gamma_j$. Finally, we assume that $L(\Gamma_j) \leq 4\pi$ and $\kappa(p_j)=1$ for some point $p_j \in \Gamma_j$. Then $\sup_{\Gamma_j} \kappa \to 1$ as $j \to \infty$.
\end{lemma}

\textbf{Proof.} 
Suppose that there exists a real number $a>0$ such that $\sup_{\Gamma_j} \kappa \geq 1+2a$ for $j$ large. We can find a segment $\tilde{\Gamma}_j \subset \Gamma_j$ such that the geodesic curvature increases from $1+a$ to $1+2a$ along $\tilde{\Gamma}_j$. Using our assumptions, we obtain $\limsup_{j \to \infty} \sup_{\tilde{\Gamma}_j} |\frac{d\kappa}{ds}| < \infty$ and $\limsup_{j \to \infty} \sup_{\tilde{\Gamma}_j} |\frac{d^2\kappa}{ds^2}| < \infty$. Since $\kappa$ varies between $1+a$ and $1+2a$ along $\tilde{\Gamma}_j$, we must have $\liminf_{j \to \infty} L(\tilde{\Gamma}_j) > 0$. On the other hand, we have $L(\tilde{\Gamma}_j) \leq L(\Gamma_j) \leq 4\pi$. Hence, after passing to a subsequence, the curves $\tilde{\Gamma}_j$ converge in $C^3$ to a curve $\hat{\Gamma}$. The geodesic curvature of the limiting curve $\hat{\Gamma}$ increases from $1+a$ to $1+2a$ as we travel along the curve $\hat{\Gamma}$. Finally, at each point on $\hat{\Gamma}$, the inscribed radius is at least $\frac{1}{\hat{\kappa}}$, where $\hat{\kappa}$ denotes the geodesic curvature of $\hat{\Gamma}$. By Lemma \ref{const.curvature}, $\hat{\kappa}$ is constant. This contradicts the fact that $\hat{\kappa}$ varies between $1+a$ and $1+2a$. \\

\begin{lemma} 
\label{inf}
Let $\Gamma_j$ be a sequence of (possibly non-closed) embedded curves in the plane with the property that $\kappa > 0$, $|\frac{d\kappa}{ds}| \leq C_\# \, (\kappa+2\Theta)^2$, and $|\frac{d^2\kappa}{ds^2}| \leq C_\# \, (\kappa+2\Theta)^3$. Moreover, suppose that the inscribed radius is at least $\frac{1}{(1+\frac{1}{j}) \, \kappa}$ at each point on $\Gamma_j$. Finally, we assume that $L(\Gamma_j) \leq 4\pi$ and $\kappa(p_j)=1$ for some point $p_j \in \Gamma_j$. Then $\inf_{\Gamma_j} \kappa \to 1$ as $j \to \infty$.
\end{lemma}

\textbf{Proof.} 
Suppose that there exists a real number $a>0$ such that $\inf_{\Gamma_j} \kappa \leq 1-2a$ for $j$ large. We can find a segment $\tilde{\Gamma}_j \subset \Gamma_j$ such that the geodesic curvature decreases from $1-a$ to $1-2a$ along $\tilde{\Gamma}_j$. Using our assumptions, we obtain $\limsup_{j \to \infty} \sup_{\tilde{\Gamma}_j} |\frac{d\kappa}{ds}| < \infty$ and $\limsup_{j \to \infty} \sup_{\tilde{\Gamma}_j} |\frac{d^2\kappa}{ds^2}| < \infty$. Since $\kappa$ varies between $1-a$ and $1-2a$ along $\tilde{\Gamma}_j$, we must have $\liminf_{j \to \infty} L(\tilde{\Gamma}_j) > 0$. On the other hand, we have $L(\tilde{\Gamma}_j) \leq L(\Gamma_j) \leq 4\pi$. Hence, after passing to a subsequence, the curves $\tilde{\Gamma}_j$ converge in $C^3$ to a curve $\hat{\Gamma}$. The geodesic curvature of the limiting curve $\hat{\Gamma}$ decreases from $1-a$ and $1-2a$ as we travel along the curve $\hat{\Gamma}$. Finally, at each point on $\hat{\Gamma}$, the inscribed radius is at least $\frac{1}{\hat{\kappa}}$, where $\hat{\kappa}$ denotes the geodesic curvature of $\hat{\Gamma}$. By Lemma \ref{const.curvature}, $\hat{\kappa}$ is constant. This contradicts the fact that $\hat{\kappa}$ varies between $1-a$ and $1-2a$. \\

\begin{proposition} 
\label{closing.up}
Let $\Gamma_j$ be a sequence of (possibly non-closed) embedded curves in the plane with the property that $\kappa > 0$, $|\frac{d\kappa}{ds}| \leq C_\# \, (\kappa+2\Theta)^2$, and $|\frac{d^2\kappa}{ds^2}| \leq C_\# \, (\kappa+2\Theta)^3$. Moreover, suppose that the inscribed radius is at least $\frac{1}{(1+\frac{1}{j}) \, \kappa}$ at each point on $\Gamma_j$, and the outer radius is at least $\frac{\alpha}{\kappa}$ at each point on $\Gamma_j$. Finally, we assume that $\kappa(p_j)=1$ for some point $p_j \in \Gamma_j$. Then $L(\Gamma_j) < 3\pi$ for $j$ large, and we have $\lim_{j \to \infty} \sup_{\Gamma_j} |\kappa-1| = 0$.
\end{proposition}

\textbf{Proof.} 
We first show that $L(\Gamma_j) < 3\pi$ for $j$ large. Suppose by contradiction that $L(\Gamma_j) \geq 3\pi$ for all $j$. By shortening $\Gamma_j$ if necessary, we can arrange that $L(\Gamma_j) = 3\pi$ for all $j$. Let $\gamma_j: [0,3\pi] \to \mathbb{R}^2$ be a parametrization of $\Gamma_j$ by arclength. It follows from Lemma \ref{sup} and Lemma \ref{inf} that the geodesic curvature of $\Gamma_j$ is close to $1$ when $j$ is sufficiently large. This implies that $\gamma_j(2\pi) - \gamma_j(0) \to 0$ as $j \to \infty$. Let us pick a sequence of numbers $s_j \in [0,3\pi]$ such that $s_j \to 2\pi$ as $j \to \infty$ and the function $s \mapsto |\gamma_j(s) - \gamma_j(0)|^2$ has a local minimum at $s_j$. Then the vector $\gamma_j(s_j)-\gamma_j(0)$ is parallel to $\nu_j(s_j)$. Consequently, we have 
\[|\langle \gamma_j(s_j)-\gamma_j(0),\nu_j(s_j) \rangle| = |\gamma_j(s_j)-\gamma_j(0)|.\] 
On the other hand, we know that the inscribed radius and the outer radius of $\Gamma_j$ are at least $\frac{\alpha}{\kappa}$. This implies 
\[\frac{1}{2} \, \kappa_j(s_j) \, |\gamma_j(s_j)-\gamma_j(0)|^2 \geq \alpha \, |\langle \gamma_j(s_j)-\gamma_j(0),\nu_j(s_j) \rangle|.\] 
Putting these facts together, we obtain 
\[\frac{1}{2} \, \kappa_j(s_j) \, |\gamma_j(s_j)-\gamma_j(0)| \geq \alpha.\] 
But $\kappa_j(s_j) \to 1$ and $|\gamma_j(s_j)-\gamma_j(0)| \to 0$ as $j \to \infty$, so we arrive at a contradiction. Consequently, we must have $L(\Gamma_j) < 3\pi$ when $j$ is sufficiently large. Using Lemma \ref{sup} and Lemma \ref{inf}, we obtain $\lim_{j \to \infty} \sup_{\Gamma_j} \kappa = 1$ and $\lim_{j \to \infty} \inf_{\Gamma_j} \kappa = 1$. This completes the proof. \\

\begin{corollary}
\label{closing.up.2}
We can find a number $\delta>0$ with the following property: Suppose that $\Gamma$ is a 
(possibly non-closed) embedded curve in the plane with the property that $\kappa > 0$, $|\frac{d\kappa}{ds}| \leq C_\# \, (\kappa+2\Theta)^2$, and $|\frac{d^2\kappa}{ds^2}| \leq C_\# \, (\kappa+2\Theta)^3$. Moreover, suppose that the inscribed radius is at least $\frac{1}{(1+\delta) \, \kappa}$ at each point on $\Gamma$, and the outer radius is at least $\frac{\alpha}{\kappa}$ at each point on $\Gamma$. Finally, we assume that $\kappa=1$ at some point $p \in \Gamma$. Then $L(\Gamma) < 3\pi$ and $\sup_\Gamma |\kappa-1| \leq \frac{1}{100}$. 
\end{corollary}

Note that the constant $\delta$ will depend only on the constants $\alpha$ and $C_\#$, which have already been chosen.

In the following, we define $\theta_0 = 10^{-6} \, \min\{\alpha,\frac{1}{C_\# \, \Theta^3}\}$.

\begin{proposition}
\label{smoothing}
We can choose $\delta$ small enough so that the following holds: Consider a family of simple closed, convex curves $\Gamma_t$, $t \in (-2\theta_0,0]$, in the plane which evolve by curve shortening flow. Assume that, for each $t \in (-2\theta_0,0]$, the curve $\Gamma_t$ satisfies the derivative estimates $|\frac{d\kappa}{ds}| \leq C_\# \, (\kappa+2\Theta)^2$ and $|\frac{d^2\kappa}{ds^2}| \leq C_\# \, (\kappa+2\Theta)^3$. Moreover, we assume that the inscribed radius is at least $\frac{1}{(1+\delta) \, \kappa}$ at each point on $\Gamma_t$, and the outer radius is at least $\frac{\alpha}{\kappa}$ at each point on $\Gamma_t$. Finally, we assume that the geodesic curvature of $\Gamma_0$ is equal to $1$ somewhere. Then the curve $\Gamma_0$ satisfies $\sum_{l=1}^{18} |\nabla^l \kappa| \leq \frac{1}{1000}$. Moreover, we have $\sup_{\Gamma_{-\theta_0}} \kappa \leq 1-\frac{\theta_0}{4}$.
\end{proposition}

\textbf{Proof.} 
Suppose that the assertion is false, and consider a sequence of counterexamples. These counterexamples converge to a smooth solution of the curve shortening flow which is defined for $t \in (-2\theta_0,0]$. The limiting solution is a family of homothetically shrinking circles. This gives a contradiction. \\

Proposition \ref{choice.of.delta} follows by combining Corollary \ref{closing.up} and Proposition \ref{smoothing}. 

\section{Proof of Proposition \ref{choice.of.hat.delta}} 

We again argue by contradiction. Let us fix $\theta_0$ and $\delta$ as above, and suppose that there is no real number $\hat{\delta} \in (0,\delta)$ for which the conclusion of Proposition \ref{choice.of.delta} holds. By taking a sequence of counterexamples and passing to the limit, we obtain a smooth solution $\Gamma_t$, $t \in (-2\theta_0,0]$, to the curve shortening flow with the property that $\sup_{\Gamma_t} \frac{\mu}{\kappa} \leq 1+\delta$ for each $t \in (-2\theta_0,0]$ and $\sup_{\Gamma_0} \frac{\mu}{\kappa} = 1+\delta$. (As usual, $\mu$ denotes the reciprocal of the inscribed radius and $\kappa$ denotes the geodesic curvature.) The geodesic curvature satisfies the evolution equation 
\[\frac{\partial}{\partial t} \kappa = \Delta \kappa + \kappa^3.\] 
Moreover, $\mu$ satisfies the inequality 
\[\frac{\partial}{\partial t} \mu \leq \Delta \mu + \kappa^2 \, \mu - \frac{2}{\mu-\kappa} \, |\nabla \mu|^2\] 
on the set $\{\mu>\kappa\}$, where $\Delta \mu$ is interpreted in the sense of distributions (see \cite{Brendle2}, Proposition 2.3). In particular, the function $(1+\delta) \, \kappa - \mu$ is nonnegative and satisfies the inequality 
\[\frac{\partial}{\partial t} ((1+\delta) \, \kappa - \mu) \geq \Delta ((1+\delta) \, \kappa - \mu) + \kappa^2 \, ((1+\delta) \, \kappa - \mu) + \frac{2}{\mu-\kappa} \, |\nabla \mu|^2\] 
on the set $\{\mu>\kappa\}$. Since $\inf_{\Gamma_0} ((1+\delta) \, \kappa - \mu) = 0$, the function $(1+\delta) \, \kappa - \mu$ vanishes identically by the strict maximum principle. This implies $\nabla \mu = 0$, hence $\nabla \kappa = 0$. Therefore, our solution is a family of shrinking circles. In that case, we have $\mu = \kappa$, which contradicts the fact that $(1+\delta) \, \kappa - \mu = 0$. 

\section{Proof of Proposition \ref{cylindrical.estimate}}

Let $\delta$ and $\hat{\delta}$ be chosen such that Theorem \ref{choice.of.delta} holds. In the following, we put 
\[f_{\delta,\sigma} = H^{\sigma-1} \, (\mu-(1+\delta) \, H) - C_1(\delta/2),\] 
where $\mu$ denotes the reciprocal of the inscribed radius and $C_1(\delta)$ is the constant in the convexity estimate of Huisken and Sinestrari (see Proposition \ref{huisken.sinestrari.convexity} above). The following result was established in \cite{Brendle2}:

\begin{proposition}[cf. \cite{Brendle2}]
\label{Lp.bound.for.f}
We can find a constant $c_0$, depending only on $\delta$ and the initial data, with the following property: if $p \geq \frac{1}{c_0}$ and $\sigma \leq c_0 \, p^{-\frac{1}{2}}$, then we have 
\[\frac{d}{dt} \bigg ( \int_{M_t} f_{\delta,\sigma,+}^p \bigg ) \leq C \, \sigma \, p \int_{M_t} f_{\delta,\sigma,+}^p + \sigma \, p \, K_0^p \int_{M_t} |A|^2\] 
except if $t$ is a surgery time. Here, $C$ and $K_0$ depend only on $\delta$ and the initial data, but not on $\sigma$ and $p$.
\end{proposition}

We next analyze the behavior of $f_{\delta,\sigma}$ under surgery. 

\begin{lemma} 
\label{behavior.under.surgery}
The integral $\int_{M_t} f_{\delta,\sigma,+}^p$ does not increase under surgery. 
\end{lemma}

\textbf{Proof.} 
Consider a surgery time $t_0$. By assumption, each surgery is being performed on an $(\hat{\alpha},\hat{\delta},\varepsilon,L)$-neck with $\varepsilon \leq \bar{\varepsilon}$ and $L \geq 1000 \, \Lambda$. Hence, Theorem \ref{properties.of.surgery} implies that the inscribed radius of $M_{t_0+}$ is at least $\frac{1}{(1+\delta) \, H}$ in the surgically modified region. In other words, we have $f_{\delta,\sigma} \leq 0$ in the surgically modified region of $M_{t_0+}$. Consequently, we have $\int_{M_{t_0+}} f_{\delta,\sigma,+}^p \leq \int_{M_{t_0-}} f_{\delta,\sigma,+}^p$, as claimed. \\

Combining Proposition \ref{Lp.bound.for.f} and Lemma \ref{behavior.under.surgery}, we can draw the following conclusion: 

\begin{proposition}
We can find a constant $c_0$, depending only on $\delta$ and the initial data, with the following property: if $p \geq \frac{1}{c_0}$ and $\sigma \leq c_0 \, p^{-\frac{1}{2}}$, then we have 
\[\int_{M_t} f_{\delta,\sigma,+}^p \leq C\] 
for all $t$. Here, $C$ is a constant that depends on $\delta$, $\sigma$, $p$, and the initial data.
\end{proposition}

We can now use Stampacchia iteration to show that $f_{\delta,\sigma} \leq C$, where $\sigma$ and $C$ depend only on $\delta$ and the initial data. This completes the proof of Proposition \ref{cylindrical.estimate}.

\section{Proof of the Neck Detection Lemma (Version A)}

\label{proof.of.neck.detection.lemma.version.a}

The proof is by contradiction. Suppose that the assertion is false. Then there exists a sequence of flows $\mathcal{M}_j$ and a sequence of points $(p_j,t_j)$ with the following properties: 
\begin{itemize} 
\item For each $j$, $\mathcal{M}_j$ is a mean curvature flow with surgery satisfying Assumption \ref{a.priori.assumptions}. 
\item $H_{\mathcal{M}_j}(p_j,t_j) \geq \max \{j,\frac{H_{1,j}}{\Theta}\}$ and $\frac{\lambda_{1,\mathcal{M}_j}(p_j,t_j)}{H_{\mathcal{M}_j}(p_j,t_j)} \leq \frac{1}{j}$. 
\item The neighborhood $\hat{\mathcal{P}}_{\mathcal{M}_j}(p_j,t_j,L_0+4,2\theta_0)$ does not contain surgeries. 
\item The point $p_j$ does not lie at the center of an $(\hat{\alpha},\hat{\delta},\varepsilon_0,L_0)$-neck in the surface $M_{t_j,j}$.
\end{itemize}
For each $j$, we put 
\begin{align*} 
\rho_j 
&= \min \Big \{ \inf \{d_{g(t_j)}(p_j,x) \, H_{\mathcal{M}_j}(p_j,t_j): \\ 
&\hspace{30mm} x \in M_{t_j,j}, \, H_{\mathcal{M}_j}(x,t_j) > 4 \, H_{\mathcal{M}_j}(p_j,t_j)\},L_0+2 \Big \}. 
\end{align*}
Using Proposition \ref{gradient.estimate}, we obtain  
\[\liminf_{j \to \infty} \rho_j > 0.\] 
By definition of $\rho_j$, we have 
\[H_{\mathcal{M}_j}(x,t_j) \leq 4 \, H_{\mathcal{M}_j}(p_j,t_j)\] 
for all points $x \in M_{t_j,j}$ satisfying $d_{g(t_j)}(p_j,x) < \rho_j \, H_{\mathcal{M}_j}(p_j,t_j)^{-1}$. Using Proposition \ref{gradient.estimate}, we obtain 
\[H_{\mathcal{M}_j}(x,t) \leq 8 \, H_{\mathcal{M}_j}(p_j,t_j)\] 
for all points $(x,t) \in \hat{\mathcal{P}}_{\mathcal{M}_j}(p_j,t_j,\rho_j,2\theta_0)$.

We next consider the restriction of the flow $\mathcal{M}_j$ to the parabolic region $\hat{\mathcal{P}}_{\mathcal{M}_j}(p_j,t_j,\rho_j,2\theta_0)$. Let us shift $(p_j,t_j)$ to $(0,0)$ and dilate the surface by the factor $H_{\mathcal{M}_j}(p_j,t_j)$. As a result, we obtain a flow $\tilde{\mathcal{M}_j}$ which is defined in the parabolic region $\mathcal{P}_{\tilde{\mathcal{M}}_j}(0,0,\rho_j,2\theta_0)$ and satisfies $H_{\tilde{M}_j}(0,0)=1$ and $\lambda_{1,\tilde{\mathcal{M}}_j}(0,0) \leq \frac{1}{j}$. Furthermore, the mean curvature of $\tilde{\mathcal{M}}_j$ is at most $8$ everywhere in the parabolic region $\mathcal{P}_{\tilde{\mathcal{M}}_j}(0,0,\rho_j,2\theta_0)$.

After passing to a subsequence, the flows $\tilde{\mathcal{M}_j}$ converge smoothly to a limit flow $\hat{\mathcal{M}}$. The limit flow is defined in a parabolic region $\mathcal{P}(0,0,\rho,2\theta_0)$, where $\rho = \lim_{j \to \infty} \rho_j > 0$. Moreover, the limit flow $\hat{\mathcal{M}}$ satisfies $H(0,0)=1$ and $\lambda_1(0,0) = 0$. Finally, the mean curvature of $\hat{\mathcal{M}}$ is at most $8$ everywhere in the parabolic region $\mathcal{P}(0,0,\rho,2\theta_0)$.

By the strict maximum principle, the limit flow $\hat{\mathcal{M}}$ splits as a product. In other words, we can find a one-parameter family of curves $\Gamma_t$, $t \in (-2\theta_0,0]$, such that $\hat{M}_t \subset \Gamma_t \times \mathbb{R}$. We may assume that the curve $\Gamma_t$ coincides with the image of $\hat{M}_t$ under the projection from $\mathbb{R}^3$ to $\mathbb{R}^2$. Note that the curves $\Gamma_t$ need not be closed. 

It follows from the Huisken-Sinestrari convexity estimate that the second fundamental form of the limiting solution $\hat{\mathcal{M}}$ is nonnegative. Hence, the curve $\Gamma_t$ has positive geodesic curvature. Since the original flow $\mathcal{M}_j$ satisfies the gradient estimate in Proposition \ref{gradient.estimate}, the curve $\Gamma_t$ satisfies the derivative estimates $|\frac{d\kappa}{ds}| \leq C_\# \, (\kappa+\Theta)^2$ and $|\frac{d^2\kappa}{ds^2}| \leq C_\# \, (\kappa+\Theta)^3$. Furthermore, since the original flow $\mathcal{M}_j$ satisfies the cylindrical estimate in Proposition \ref{cylindrical.estimate}, the limit flow $\hat{M}_t$ is $\frac{1}{1+\delta}$-noncollapsed. Hence, the inscribed radius of $\Gamma_t$ is at least $\frac{1}{(1+\delta) \, \kappa}$, where $\kappa$ denotes the geodesic curvature of $\Gamma_t$. Furthermore, the outer radius of $\Gamma_t$ is at least $\frac{\alpha}{\kappa}$ at each point on $\Gamma_t$. 

Finally, the curve $\Gamma_0$ passes through the origin, and the geodesic curvature of $\Gamma_0$ is equal to $1$ at the origin. Hence, Proposition \ref{choice.of.delta} implies that $\Gamma_0$ has length at most $3\pi$, and $\sup_{\Gamma_0} |\kappa-1| \leq \frac{1}{100}$. Moreover, Proposition \ref{gradient.estimate} implies that each curve $\Gamma_t$ contains a point where the geodesic curvature is between $\frac{1}{2}$ and $2$. Applying Proposition \ref{choice.of.delta} to a scaled copy of $\Gamma_t$, we conclude that each curve $\Gamma_t$ has length at most $6\pi$. 

At this point, we distinguish two cases: 

\textit{Case 1:} Suppose first that $0 < \rho < L_0+2$. Then $\rho_j < L_0+2$ for $j$ large. From this, we deduce that $\sup_{\tilde{M}_{0,j}} H \geq 4$ for $j$ large. Using this fact and the gradient estimate, we obtain $\sup_{\hat{M}_0} H \geq 2$. Consequently, $\sup_{\Gamma_0} \kappa \geq 2$, where $\kappa$ denotes the geodesic curvature of $\Gamma_0$. On the other hand, we have established earlier that $\sup_{\Gamma_0} |\kappa-1| \leq \frac{1}{100}$. This is a contradiction. 

\textit{Case 2:} We now assume that $\rho = L_0+2 > 100$. Since $\Gamma_t$ has length at most $6\pi$, we conclude that $\Gamma_t$ must be a closed curve. Hence, the curves $\Gamma_t$ are simple closed, convex curves in the plane, which evolve by curve shortening flow. 

By Proposition \ref{choice.of.delta}, the curve $\Gamma_0$ satisfies $\sum_{l=1}^{18} |\nabla^l \kappa| \leq \frac{1}{1000}$. Moreover, we have $\sup_{\Gamma_{-\theta_0}} \kappa \leq 1-\frac{\theta_0}{4}$. Finally, Proposition \ref{choice.of.hat.delta} implies that, for each point on $\Gamma_0$, the inscribed radius is at least $\frac{1}{(1+\hat{\delta}) \, \kappa}$. 

If $j$ is sufficiently large, we can find a region $N_j \subset \{x \in M_{t_j,j}: d_{g(t_j)}(p_j,t_j) \leq (L_0+1) \, H(p_j,t_j)^{-1}\}$ such that $\text{\rm dist}_{C^{20}}(H(p_j,t_j) \, (N_j-p_j),\Gamma_0 \times [-L_0,L_0]) < \varepsilon_0$. We again divide the discussion into two subcases:

\textit{Subcase 2.1:} Suppose that, for $j$ large, the region 
\[\{x + a \, \nu(x): x \in M_{t_j,j}, \, d_{g(t_j)}(p_j,x) \leq (L_0+1) \, H(p_j,t_j)^{-1}, \, a \in (0,2\hat{\alpha} \, H(p_j,t_j)^{-1})\}\] 
is disjoint from $M_{t_j,j}$. Consequently, the point $p_j$ lies at the center of an $(\hat{\alpha},\hat{\delta},\varepsilon_0,L_0)$-neck in $M_{t_j,j}$ if $j$ is sufficiently large. This contradicts our assumption.

\textit{Subcase 2.2:} Suppose that, for $j$ large, the region 
\[\{x + a \, \nu(x): x \in M_{t_j,j}, \, d_{g(t_j)}(p_j,x) \leq (L_0+1) \, H(p_j,t_j)^{-1}, \, a \in (0,2\hat{\alpha} \, H(p_j,t_j)^{-1})\}\] 
does intersect $M_{t_j,j}$. In this case, we can find a sequence of points $x_j \in M_{t_j,j}$ and a sequence of numbers $a_j \in (0,2\hat{\alpha} \, H(p_j,t_j)^{-1})$ such that $d_{g(t_j)}(p_j,x_j) \leq (L_0+1) \, H(p_j,t_j)^{-1}$ and $z_j := x_j + a_j \, \nu(x_j) \in M_{t_j,j}$. We next observe that $H(x_j,t_j) \leq 4 \, H(p_j,t_j)$ by definition of $\rho_j$. Hence, the outer radius of the surface $M_{t_j,j}$ at the point $x_j$ is at least $\alpha \, H(x_j,t_j)^{-1} \geq \frac{\alpha}{4} \, H(p_j,t_j)^{-1}$. Consequently, $a_j \geq \frac{\alpha}{2} \, H(p_j,t_j)^{-1}$. 

We now let $\tau_j = t_j - \theta_0 \, H(p_j,t_j)^{-2}$. If $j$ is sufficiently large, we may write $z_j = y_j + b_j \, \nu(y_j)$, where $(y_j,\tau_j) \in \hat{\mathcal{P}}_{\mathcal{M}_j}(p_j,t_j,L_0+1,2\theta_0)$ and $0 \leq b_j \leq a_j < 2\hat{\alpha} \, H(p_j,t_j)^{-1}$. On the other hand, since $\sup_{\Gamma_{-\theta_0}} \kappa \leq 1-\frac{\theta_0}{4}$, we have $H(y_j,\tau_j) \leq (1-\frac{\theta_0}{8}) \, H(p_j,t_j)$. This implies that the outer radius of the surface $M_{\tau_j,j}$ at the point $y_j$ is at least 
\[\alpha \, H(y_j,\tau_j)^{-1} \geq \frac{\alpha}{1-\frac{\theta_0}{8}} \, H(p_j,t_j)^{-1} = \hat{\alpha} \, H(p_j,t_j)^{-1} > \frac{b_j}{2}.\] 
Consequently, the point $z_j = y_j + b_j \, \nu(y_j)$ does not lie in the region enclosed by $M_{\tau_j,j}$. This contradicts the fact that $z_j \in M_{t_j,j}$. This completes the proof.

\section{Proof of Proposition \ref{7.12}}

\label{proof.of.7.12}

By assumption, the parabolic neighborhood $\hat{\mathcal{P}}(p_1,t_1,\tilde{L}+4,2\theta_0)$ contains a point which belongs to a surgery region. Consequently, we can find a surgery time $t_0 \in [t_1-2\theta_0 \, H(p_1,t_1)^{-2},t_1)$ and a point $q_1 \in M_{t_1}$ such that the following holds: 
\begin{itemize}
\item $d_{g(t_1)}(p_1,q_1) \leq (\tilde{L}+4) \, H(p_1,t_1)^{-1}$.
\item If we follow the point $q_1 \in M_{t_1}$ back in time, then the corresponding point $q_0 \in M_{t_0+}$ lies in the region modified by surgery at time $t_0$.
\end{itemize}
Let us consider the region modified by surgery at time $t_0$, and let $U_0$ denote the connected component of this set that contains the point $q_0$. In other words, $U_0 \subset M_{t_0+}$ is a cap that was inserted at time $t_0$. We next define $V_0 = \{x \in M_{t_0+}: \text{\rm dist}_{g(t_0+)}(U_0,x) \leq 1000 \, H_1^{-1}\}$. Clearly, $V_0$ is diffeomorphic to a disk. Let 
\[D = \{y \in \mathbb{R}^3: \text{\rm there exists a point $x \in V_0$ such that $|y - x| < \frac{\alpha}{1000} \, H_1^{-1}$}\}.\] 
Arguing as in Proposition \ref{separation.of.surgery.regions} above, we can show that, for every surgery time $t>t_0$, the set $D$ is disjoint from the region modified by surgery at time $t$. Consequently, the surfaces $M_t \cap D$ form a regular mean curvature flow for $t > t_0$. In other words, the surfaces $M_t \cap D$ evolve smoothly for $t > t_0$, but we allow the possibility that some components of $M_t \cap D$ may disappear as a result of surgeries in other regions.

At each point on $V_0 \subset M_{t_0+}$, the mean curvature is at most $20 \, H_1$. We now follow the surface $V_0 \subset M_{t_0+}$ forward in time. This gives a one-parameter family of surfaces which are all diffeomorphic to a disk. It follows from Proposition \ref{gradient.estimate} that, for $t \in (t_0,t_0+2\theta_0 \, H_1^{-2}]$, the resulting surfaces remain inside the region $D$ and have mean curvature at most $40 \, H_1$. Moreover, since $q_1 \in M_{t_1}$, the resulting surfaces cannot disappear before time $t_1$.

Let $V_1 \subset M_{t_1}$ denote the region in $M_{t_1}$ which is obtained by following the region $V_0 \subset M_{t_0+}$ forward in time. Clearly, $V_1$ is diffeomorphic to a disk, and the mean curvature is at most $40 \, H_1$ at each point in $V_1$. Since $q_0 \in V_0$, we have $q_1 \in V_1$. Furthermore, since $\text{\rm dist}_{g(t_0+)}(q_0,\partial V_0) \geq 1000 \, H_1^{-1}$, we obtain $\text{\rm dist}_{g(t_1)}(q_1,\partial V_1) \geq 500 \, H_1^{-1}$. From this, we deduce that 
\[\{x \in M_{t_1}: d_{g(t_1)}(q_1,x) \leq 500 \, H_1^{-1}\} \subset V_1.\] 
Hence, if we put $V := V_1$, then $V$ has the required properties.

\section{Proof of Proposition \ref{7.19}}

\label{proof.of.7.19}

To fix notation, let $U \subset \{x \in \mathbb{R}^3: y_1 \leq x_3 \leq y_2\}$ denote the region enclosed by $\Sigma$. Moreover, let $\nu$ denote the outward-pointing unit normal to $U$. 

Suppose by contradiction that there exists a point $\bar{x} \in \Sigma$ such that $\bar{x}_3 \in [y_1+1,y_2]$, $\langle \nu(\bar{x}),e_3 \rangle \geq 0$, and $H(\bar{x}) \leq \frac{\alpha}{100}$. The noncollapsing assumption implies that there exists a ball $B \subset \mathbb{R}^3$ of radius $100$ such that 
\[B \cap \{x \in \mathbb{R}^3: y_1 \leq x_3 \leq y_2\} \subset U.\] 
This implies 
\[B \cap \{x \in \mathbb{R}^3: x_3=\bar{x}_3-1\} \subset U \cap \{x \in \mathbb{R}^3: x_3=\bar{x}_3-1\}.\] 
Since $\langle \nu(\bar{x}),e_3 \rangle \geq 0$, the set 
\[B \cap \{x \in \mathbb{R}^3: x_3=\bar{x}_3-1\}\] 
is a disk of radius at least $\sqrt{100^2-99^2} > 10$. On the other hand, our assumptions imply that the set $U \cap \{x \in \mathbb{R}^3: x_3=\bar{x}_3-1\}$ is contained in a disk of radius $10$. This is a contradiction.

\end{document}